\documentclass{article}
\usepackage{ijcai25}

\usepackage{times}
\usepackage{soul}
\usepackage{url}
\usepackage[hidelinks]{hyperref}
\usepackage[utf8]{inputenc}
\usepackage[small]{caption}
\usepackage{graphicx}
\usepackage{amsmath}
\usepackage{amsthm}
\usepackage{booktabs}
\usepackage{algorithm}
\usepackage{algorithmic}
\usepackage[switch]{lineno}

\usepackage{bm}
\usepackage{amsfonts}
\usepackage{amssymb}
\usepackage[capitalise]{cleveref}

\newtheorem{condition}{Condition}

\newtheorem{remark}{Remark}

\usepackage{subfig}

\newtheorem{theorem}{Theorem}

\title{A Priori Estimation of the Approximation, Optimization and Generalization Errors of Random Neural Networks for Solving Partial Differential Equations}

\author{
Xianliang Xu$^4$
\and
Ye Li$^{1,2,3*}$\and
Zhongyi Huang$^4$\thanks{Corresponding author: Ye Li $<$yeli20@nuaa.edu.cn$>$ and Zhongyi Huang $<$zhongyih@tsinghua.edu.cn$>$}\\
\affiliations
$^1$College of Computer Science and Technology, Nanjing University of Aeronautics and Astronautics\\
$^2$MIIT Key Laboratory of Pattern Analysis and Machine Intelligence, Nanjing\\
$^3$State Key Laboratory for Novel Software Technology, Nanjing University\\
$^4$Department of Mathematical Sciences, Tsinghua University\\
\emails
xuxl19@mails.tsinghua.edu.cn,
yeli20@nuaa.edu.cn,
zhongyih@tsinghua.edu.cn
}

\begin{document}

\maketitle

\begin{abstract}
In recent years, neural networks have achieved remarkable progress in various fields and have also drawn much attention in applying them on scientific problems. A line of methods involving neural networks for solving partial differential equations (PDEs), such as Physics-Informed Neural Networks (PINNs) and the Deep Ritz Method (DRM), has emerged. Although these methods outperform classical numerical methods in certain cases, the optimization problems involving neural networks are typically non-convex and non-smooth, which can result in unsatisfactory solutions for PDEs. In contrast to deterministic neural networks, the hidden weights of random neural networks are sampled from some prior distribution and only the output weights participate in training. This makes training much simpler, but it remains unclear how to select the prior distribution. In this paper, we focus on Barron type functions and approximate them under Sobolev norms by random neural networks with clear prior distribution. In addition to the approximation error, we also derive bounds for the optimization and generalization errors of random neural networks for solving PDEs when the solutions are Barron type functions. 
\end{abstract}

\section{Introduction}
As the development of hardcore and algorithms, deep neural networks have made remarkable progress across various fields, including computer vision \cite{9}, natural language processing \cite{10}, reinforcement learning \cite{11}, and others. These successes have inspired researchers to explore the application of neural networks to scientific challenges, particularly in the modeling of physical systems. In the field of scientific computing, a long-standing problem is solving partial differential equations (PDEs) numerically, which can be hindered by the curse of dimensionality when using classical numerical methods. To tackle PDE-related problems, several neural network-based approaches have been introduced. Among these, Physics-Informed Neural Networks (PINNs) and the Deep Ritz Method (DRM) stand out. PINN incorporates prior information from PDEs into the training in solving forward and inverse problems of PDEs. Specifically, it encodes the PDE constraints into the design of the loss function, which restricts the constructed neural network to follow physical law characterized by the PDE. This framework's flexibility stems from its reliance solely on the PDE's form, making it adaptable for a wide range of PDEs. The DRM, on the other hand, incorporates the variational formulation that is an essential tool in traditional methods, into training the neural networks. This approach typically involves lower-order derivatives, potentially offering greater computational efficiency compared to PINNs. However, the DRM's utility is somewhat limited by the fact that not all PDEs possess a variational formulation. 

The success of neural networks can be partly attributed to their powerful approximation capabilities. The capacity of neural networks with a variety of activation functions to approximate different types of target functions has been extensively studied. This includes continuous functions \cite{13}, smooth functions \cite{14}, as well as Sobolev functions \cite{15,16,17,18} and Barron functions \cite{19,20,21}, among others. In short, neural networks can approximate the aforementioned function classes with arbitrary precision.
Despite the widespread adoption and impressive approximation capabilities of neural networks in scientific computing, their practical application can encounter challenges. A fundamental issue arises in solving the optimization problems involving neural networks, which are typically non-convex and non-smooth. For example, \cite{7} demonstrates that even for simple problems of convection, reaction, and reaction-diffusion, PINN approach only works for very easy parameter regimes and fails in more challenging physical regimes. By analyzing the loss landscape of PINN, they show that the failure is not due to the limited capacity of the neural network, but rather due to the optimization difficulties caused by the PINN's soft PDE constraint. 

 Due to the limitations described before, there is a growing interest in the application of random neural networks, whose hidden weights are randomly generated and only the output weights are trainable.
Compared to deterministic neural networks, random neural networks can lead to optimization problems that can be efficiently solved. For instance, employing random neural networks in $L^2$ regression problems results in least squares problems, which possess closed-form solutions or can be addressed by various optimization algorithms. Because of the favorable properties, random neural networks have been successfully applied not only in traditional machine learning tasks but also in addressing PDE-related problems. In this work, we focus on the utilization of random neural networks in the framework of PINNs. First, we establish the approximation error of neural networks for Barron-type functions under the $H^2$ norm. Subsequently, for applying random neural networks for solving a certain class of second-order elliptic PDEs with Barron-type solutions, we provide the optimization and generalization error analysis, providing a comprehensive understanding of how random neural networks perform when applied to PDEs.

\subsection{Contributions}
The contribution of this work can be summarized as follows.
\begin{itemize}
\item We provides two approximation results for two Barron-type functions via random neural networks under the Sobolev norms. Unlike previous works that only showed the existence of the prior distributions, we give the concise forms of these distributions. 
\item When applying random neural networks for solving certain second-order elliptic PDEs whose solutions are Barron-type, we design tailored  optimization algorithms for addressing these problems. Subsequently, we perform a full error analysis that rigorously bounds the approximation, optimization, and generalization errors. In deriving the generalization bounds, we also show that projected gradient descent can achieve the faster rate $\mathcal{O}(1/n)$.

\item We have conducted numerical experiments to validate our conclusions.
\end{itemize}

\subsection{Related Works}
\paragraph{Random Neural Networks in Machine Learning.} Random neural networks, also known as Extreme Learning Machines (ELMs) in the field of machine learning, have drawn significant attention due to their special training methodology and efficiency. The concept of ELMs was first introduced in \cite{29} as a two-layer random neural network where the hidden weights are randomly assigned, and only the output weights are trained. As shown in \cite{29}, this algorithm can provide good generalization performance at extremely fast learning speed and it outperforms deterministic neural networks in regression, classification and real-world complex problems. It is also possible to extend ELM to kernel learning \cite{30}, which shows that ELM can use a wide type of feature mappings, including random hidden nodes and kernels. 

\paragraph{Random Neural Networks in Scientific Computing.} The potential of random neural networks extends beyond traditional machine learning domains, they have also been effectively utilized in scientific computing, particularly for solving problems related to PDEs.
For instance, \cite{33} has developed an efficient method based on domain decomposition and random neural networks to solve different types of PDEs, showing the first time when the neural network-based methods outperform the traditional numerical methods in low dimensions. Later, this method has been extended to high dimensions \cite{34}, producing accurate solutions to high-dimensional PDEs. 

\cite{35} has studied approximation based on single-hidden-layer feedforward and recurrent neural networks with randomly generated hidden weights under the $L^2$ norm. Similar to our work, they are also based on an integral representations of the target functions. However, the prior distribution for the hidden weights remains unclear in their study. In contrast, \cite{37} derived approximation rates and an explicit algorithm to learn a deterministic
function by a random neural network under certain Sobolev norms, but the conditions for the target functions to be approximated are challenging to verify. The work most closely related to ours is \cite{36}, which provided a full analysis of random neural networks for learning sufficiently non-degenerate Black-Scholes type models. However, their approximation results are presented under the $L^{\infty}$ norm, and their generalization analysis is grounded in methods for $L^2$ regression problems. Moreover, their methods cannot be applied to our setting, which involves unsupervised learning as opposed to the supervised problems they consider. To the best of our knowledge, this is the first article that attempts to provide theoretical support for solving PDEs within the framework of PINNs using random neural networks.

\subsection{Notations}
For $x\in \mathbb{R}^d$, $\|x\|_p$ denotes its $p$-norm ($1\leq p\leq \infty$). For the activation functions, we write $\sigma_k(x)$ for the $\text{ReLU}^k$ activation function, i.e., $\sigma_k(x):=[\max(0,x)]^k$. For given probability measure $P$ and a sequence of random variables $\{X_i\}_{i=1}^n$ distributed according to $P$, we denote the empirical measure of $P$ by $P_n$, i.e. $P_n=\frac{1}{n}\sum\limits_{i=1}^n \delta_{X_i}$.

\section{Preliminaries}
In this section, we provide some preliminaries about random neural networks and Barron functions, which are pivotal to our study. Throughout the paper, we only consider the two-layer random neural networks, i.e., a feedforward neural network with one hidden layer and randomly generated hidden weights. For brevity, we still call them random neural networks. 

To make it more precise, for $m\in \mathbb{N}$, let $W_1,\cdots, W_m$ be i.i.d. $\mathbb{R}^d$-valued random vectors and $B_1,\cdots, B_m$ be i.i.d. real-valued random variables, where $W=(W_1,\cdots, W_m) $ and $B=(B_1,\cdots, B_m)$ are independent. Then for any $\mathbb{R}^m$-valued (random) vector $A=(A_1,\cdots, A_m)$, we have the random neural network
\begin{equation}
H_{A}^{W,B}(x):=\sum\limits_{i=1}^m A_i\sigma(W_i\cdot x+B_i),
\end{equation} 
where $\sigma:\mathbb{R}\rightarrow \mathbb{R}$ is a fixed activation function. 

In addressing the $L^2$ regression with random neural networks, given training samples $\{(x_i,y_i)\}_{i=1}^n$, where $x_i$ is the input and $y_i$ is the output for $1\leq i \leq n$, we have the following objective function to be optimized.
\begin{equation}
\sum\limits_{i=1}^n \left(H_{A}^{W,B}(x_i)-y_i\right)^2.
\end{equation}
This results in a least squares problem that has a closed-form solution. In contrast, deterministic neural networks also require optimizing the weights of the hidden layers, which leads to a more complex optimization problem.

In this work, we mainly focus on the Barron-type functions, leveraging their integral representations as highlighted in \cite{38}. The Barron space with order $s>0$ is defined as
\begin{equation}
\begin{aligned}
&\mathcal{B}^s(\Omega):=\{ f:\Omega\rightarrow \mathbb{C} : \\
&\|f\|_{\mathcal{B}^s(\Omega)}:=\inf\limits_{f_{e} |\Omega=f} \int_{\mathbb{R}^d} (1+\|w\|_1)^s|\hat{f}_{e}(w)|d w<\infty ,\}
\end{aligned}
\end{equation}
where the infimum is over extensions $f_e \in L^1(\mathbb{R}^d)$ and $\hat{f}_{e}$ is the Fourier transform of ${f}_{e}$. Barron introduced this class for $s=1$ and showed that two-layer neural networks with sigmoidal activation function can achieve the approximation rate $\mathcal{O}(1/\sqrt{n})$ in the $L^2$ norm. Although the convergence rate does not suffer the curse of dimensionality, the related optimization problems are non-convex and challenging to address. Consequently, we shift our focus toward approximating Barron functions employing random neural networks. Note that we choose $1$-norm in the definition (3) just for simplicity. In the following, we assume that $\Omega$ is a subset of $[0,1]^d$.

\section{Main Results}

\subsection{Approximation Results}
Our approach to approximating Barron functions by random two-layer neural networks leverages the integral representations of these functions, as presented in \cite{1}. Specifically, for $f\in \mathcal{B}^s(\Omega)$, without loss of generality, assume that the infimum is attained at $f_e$. Then for activation function $\sigma \in L^1(\mathbb{R})$ with $\hat{\sigma}(a)\neq 0$ and some $a\neq 0$, we have
\begin{equation}
\begin{aligned}
&f(x)=  \int_{\mathbb{R}^d} e^{i\omega \cdot x}\hat{f}_e(\omega)d\omega \\
&=\int_{\mathbb{R}^d}\int_{ \mathbb{R}} \frac{1}{2\pi \hat{\sigma}(a) }\sigma(\frac{\omega}{a}\cdot x+b)\hat{f}_e(\omega)e^{-iab}dbd\omega \\
&=\int_{\mathbb{R}^d}\int_{ \mathbb{R}} \frac{1}{2\pi \hat{\sigma}(a) }\sigma(\frac{\omega}{a}\cdot x+b)|\hat{f}_e(\omega)|\cos(\theta(\omega)-ab)dbd\omega,
\end{aligned}
\end{equation}
where $\hat{f}_e(\omega)=e^{i\theta(\omega)}|\hat{f}_e(\omega)|$.

In the paper, we consider the activation function $\sigma(t)=\sum_{i=0}^{4}(-1)^{i}C_{4}^i \sigma_3(t+2-i)$, which is compactly supported on $[-2,2]$. The conclusions can be naturally extended to sigmoidal activation functions and tanh activation functions, which are the commonly used smooth activation functions for PINNs. Specific details can be found in Remark 2. The approximation results rely on the polynomial decay condition of the activation function, i.e., $|\sigma^{(k)}(t)|\leq C_p(1+|t|)^{-p}$ for $0\leq k \leq 2$ for some $p>1$. This condition also appears in \cite{1} for deriving the approximation rate for functions in $\mathcal{B}^s(\Omega)$ using two-layer deterministic neural networks.

\begin{theorem}Let $P_1$ be the uniform distribution of the domain $\{\omega\in\mathbb{R}^d: \|\omega\|_1\leq M\}$ and $P_2$ be the uniform distribution of the domain $\{b\in\mathbb{R}: |b|\leq 2M\}$ with a constant $M\geq2$. Let 
$W_1,\cdots, W_m\sim P_1$ and $B_1,\cdots, B_m\sim P_2$, then there exists a $\mathbb{R}^m$-valued, $\sigma(W,B)$-measurable vector $A=(A_1,\cdots, A_m)$ such that 
for $3<s<6$,
\begin{equation}
\begin{aligned}
&\mathbb{E}\left[ \|f(x)-f_m(x)\|_{H^2(\Omega)}^2\right]\lesssim\\
&\frac{M^{d+7-s}}{m}\|f\|_{\mathcal{B}^s(\Omega)} \|\hat{f_e}\|_{L^{\infty}(\mathbb{R}^d)}+\left(\frac{1}{M^{s-3}} \|f\|_{\mathcal{B}^s(\Omega)} \right)^2.
\end{aligned}
\end{equation}
and for $s\geq 6$,
\begin{equation}
\begin{aligned}
&\mathbb{E}\left[ \|f(x)-f_m(x)\|_{H^2(\Omega)}^2\right]\lesssim\\
&\frac{M^{d+1}}{m}\|f\|_{\mathcal{B}^s(\Omega)}\|\hat{f_e}\|_{L^{\infty}(\mathbb{R}^d)}+\left(\frac{1}{M^{s-3}} \|f\|_{\mathcal{B}^s(\Omega)} \right)^2.
\end{aligned}
\end{equation}
where 
\begin{equation}
f_m(x)=\frac{1}{m}\sum\limits_{i=1}^m A_i\sigma(W_i\cdot x+B_i)
\end{equation}
and
\begin{equation}
|A_i| \leq C  \|\hat{f_e}\|_{L^{\infty}(\mathbb{R}^d)} M^{d+1},
\end{equation}
$\lesssim $ indicates that a universal multiplicative constant is omitted.
\end{theorem} 	

\begin{remark}
By appropriately selecting the value of $M$, we can obtain that
\begin{equation}
\mathbb{E}\left[ \|f(x)-f_m(x)\|_{H^2(\Omega)}^2\right]\lesssim \left\{
\begin{aligned}
m^{-\frac{2s-6}{d+s+1}} & , & 3<s<6, \\
m^{-\frac{2s-6}{d+2s-5}} & , & s\geq 6.
\end{aligned}
\right.
\end{equation}
\end{remark}

\begin{remark}
Extending our discussion beyond the specific activation function we have chosen, we explore other common activation functions, such as the Sigmoidal function Sig(x) defined as $(1+e^{-x})^{-1}$ and the Hyperbolic tangent function Tanh(x), expressed as $\frac{e^{x}-e^{-x}}{e^{x}+e^{-x}}$. Then we can construct new activation functions $\sigma(x)=Sig(x+1)-Sig(x)$ and $\sigma(x)=Tanh(x+1)-Tanh(x)$, which satisfy that
\[|\sigma^{(k)}(t)|\leq Ce^{-|t|}\]
for $0\leq k \leq 2$, where $C$ is a universal constant.
	
With this exponential decay condition, we can deduce from the proof of Theorem 1 that
for $3<s<6$,
\begin{equation}
\begin{aligned}
&\mathbb{E}\left[ \|f(x)-f_m(x)\|_{H^2(\Omega)}^2\right]\lesssim\\
&\frac{M^{d+7-s}}{m}\|f\|_{\mathcal{B}^s(\Omega)} \|\hat{f_e}\|_{L^{\infty}(\mathbb{R}^d)}+\left(\frac{1}{M^{s-3}}+\frac{1}{e^M}  \right)^2\|f\|_{\mathcal{B}^s(\Omega)}^2.
\end{aligned}
\end{equation}
for $s\geq 6$,
\begin{equation}
\begin{aligned}
&\mathbb{E}\left[ \|f(x)-f_m(x)\|_{H^2(\Omega)}^2\right]\lesssim\\
&\frac{M^{d+1}}{m}\|f\|_{\mathcal{B}^s(\Omega)}\|\hat{f_e}\|_{L^{\infty}(\mathbb{R}^d)}+\left(\frac{1}{M^{s-3}}+\frac{1}{e^M}  \right)^2\|f\|_{\mathcal{B}^s(\Omega)}^2.
\end{aligned}
\end{equation}

Besides these activation functions, Theorem 3 in \cite{1} implies that, when $\sigma\in W^{m,\infty}(\mathbb{R})$ is a non-constant periodic function, $f$ has a similar integral representation for $\sigma$. Moreover, in this case, we do not need to truncate $b$, which may yields a better approximation rate.  
\end{remark}

\paragraph{Proof Sketch:} In the integral representations of functions within $\mathcal{B}^s(\Omega)$, the integral can be partitioned into two components, corresponding to bounded and unbounded regions, respectively. Specifically,
\begin{equation}
\begin{aligned}
&f(x)=\int_{\mathbb{R}^d} \int_{\mathbb{R}} \frac{\sigma(\omega\cdot x+b)}{2\pi \hat{\sigma}(1) }|\hat{f}(\omega)|\cos(\theta(\omega)-b)dbd\omega \\
&=\int_{\{\|\omega\|_1\leq M,|b|\leq 2M\}}  \frac{\sigma(\omega\cdot x+b)}{2\pi \hat{\sigma}(1) }|\hat{f}(\omega)|\cos(\theta(\omega)-b)dbd\omega \\
&\ +\int_{\{ \|\omega\|_1\leq M, |b|\leq 2M\}^c  } \frac{\sigma(\omega\cdot x+b)}{2\pi \hat{\sigma}(1) }|\hat{f}(\omega)|\cos(\theta(\omega)-b)dbd\omega\\
&:= f_1(x)+f_2(x), 
\end{aligned}
\end{equation}
where we write $\hat{f}$ for $\hat{f_e}$ for simplicity. For the first part $f_1(x)$, it can be approximated by a random neural network. The second part, $f_2(x)$, can be further decomposed into two terms: one corresponding to the integral over $\{\|\omega\|_1> M, b\in \mathbb{R}\}$ and another to the integral over $\{\|\omega\|_1\leq M, |b|>2M\}$. We denote these parts as $f_{21}(x)$ and $f_{22}(x)$, respectively. Given that the activation function $\sigma$ is compactly supported in $[-2,2]$, we can deduce that $\sigma(\omega\cdot x+b)=0$ in $\{\|\omega\|_1\leq M, |b|>2M\}$, since 
\[ |\omega\cdot x+b|\geq |b|-|\omega\cdot x|\geq |b|-\|\omega\|_1\|x\|_{\infty}>M\geq 2,\]
which implies $f_{22}(x)=0$.

Therefore, to ensure that $f(x)$ can be well-approximated well by the random neural network designed to approximate $f_1(x)$, it suffices to guarantee that $f_{21}(x)$ has a small enough $H^2$ norm. It is for this reason that we impose the polynomial decay condition on $\sigma$. With this condition, we can establish that
\begin{equation}
\int_{\mathbb{R}} |\sigma^{(k)}(\omega\cdot x+b)|^sdb \lesssim 1+\|\omega\|_1,
\end{equation}
which is crucial to achieve our goal.
The full proof can be found in the appendix of \cite{42}.

\begin{remark}
In \cite{36}, the authors also considered employing a random neural network to approximate the function $f$ belonging to  $\mathcal{B}^s(\Omega)$, with $W_i$ having a strictly positive Lebesgue-density $\pi_{w}$ on $\mathbb{R}^d$ and $B_i$ having a strictly positive Lebesgue-density $\pi_{b}$ on $\mathbb{R}$. Despite achieving a dimension-independent approximation rate of $\mathcal{O}(1/\sqrt{m})$, the densities $\pi_{w}$ and $\pi_{b}$ are dependent on the unknown decay of the Fourier transform of $f$, leaving the prior distributions for $W_i$ and $B_i$ unclear. For general $s$, there may not be densities $\pi_{w}$ and $\pi_{b}$ that ensure the constants in the approximation rate are finite. Nonetheless, the approximation results in \cite{36} remain valid, because the target function they consider possesses a Fourier transform that decays exponentially. Moreover, the method employed in \cite{36} is specifically tailored for the ReLU activation function and the approximation results are given in the $L^{\infty}$ norm, which renders it inapplicable to the context of PINNs.

For the $L^2$ regression problems, as shown in \cite{39}, random neural networks can achieve a rate of $\mathcal{O}(1/\sqrt{n})$ for the final prediction error, where $n$ is the sample size. This study assumes that the regression function $f$ satisfies the condition
\begin{equation}
|\hat{f}(\omega)| \leq \frac{c_1}{\|\omega\|_2^{d+1} (\log \|\omega\|_2)^2 },
\end{equation}
for all $\omega \in \mathbb{R}^d$ with $\|\omega\|_2\geq 2$ and some constant $c_1>0$. It is evident that this condition implies that 
\begin{equation}
\int_{\mathbb{R}^d} (1+\|\omega\|_2)|\hat{f}(\omega)|d\omega<\infty.
\end{equation}
Thus, the assumption made in \cite{39} is stronger than the one in our study and may be more challenging to verify.
\end{remark}	

Clearly, the approximation rates presented in Theorem 1 are affected by the curse of dimensionality, which is attributed to the unknown decay of the Fourier transform of the target function. When the smoothness parameter $s$ is sufficiently large, we can choose not to truncate $\omega$ and $b$, which may lead to better results. Specifically, assume that $s>\alpha+\beta+5$ with $\alpha>d,\beta>1$, then we can derive the following representation for $f\in \mathcal{B}^s(\Omega)$.
\begin{equation}
\begin{aligned}
&f(x)=  \int_{\mathbb{R}^d}\int_{\mathbb{R}} \frac{\sigma(\omega\cdot x+b)}{2\pi \hat{\sigma}(1) }|\hat{f}(\omega)|\cos(\theta(\omega)-b)dbd\omega\\
&= \int_{\mathbb{R}^d}\int_{\mathbb{R}} \frac{\sigma(\omega\cdot x+b)}{2\pi \hat{\sigma}(1) } |\hat{f}(\omega)|\cos(\theta(\omega)-b) \\ 
&\quad\cdot\frac{(1+\|\omega\|_2)^{\alpha}(1+|b|)^{\beta}}{C_{\alpha}C_{\beta}} \frac{C_{\alpha}C_{\beta}}{(1+\|\omega\|_2)^{\alpha} (1+|b|)^{\beta}}dbd\omega\\
&=\mathbb{E}_{P_1(\omega), P_2(b)}\left[ \frac{1}{2\pi \hat{\sigma}(1) }\frac{\sigma(\omega\cdot x+b)|\hat{f}(\omega)|\cos(\theta(\omega)-b)}{p_1(\omega)p_2(b)}\right] ,
\end{aligned}
\end{equation}
where $P_1(\omega)=p_1(\omega)d\omega, P_2(b)=p_2(b)db$ are probability measures with density functions $p_1(\omega)$ and $p_2(b)$ in $\mathbb{R}^d$ and $\mathbb{R}$, respectively, defined as $p_1(\omega):=\frac{C_{\alpha}}{(1+\|\omega\|_2)^{\alpha}}$ and $p_2(b):=\frac{C_{\beta}}{(1+|b|)^{\beta}}$.

Building upon our previous discussion, we now consider a scenario where the target function $f$ belongs to a smaller function space. Specifically, we assume that $f$ satisfies the following condition.

\begin{condition}
Given function $f$, there is a function $f_e\in L^1(\mathbb{R}^d)$ such that $f_e|_{\Omega}=f$ and 
\begin{equation}
\int_{\mathbb{R}^d} (1+\|\omega\|_2)^s |\hat{f}_e(\omega)|^2d\omega <\infty
\end{equation}
with $s>\alpha+\beta+5$ and $\alpha >d+4,1<\beta<5$.
\end{condition}

Given the integral representation of the function $f$ that satisfies Condition 1, i.e., (17), we can derive the following approximation results for $f$.

\begin{theorem}
Let $W_1,\cdots,W_m\sim P_1$ and $B_1,\cdots, B_m\sim P_2$, then there exists a $\mathbb{R}^m$-valued, $\sigma(W,B)$-measurable vector $A=(A_1,\cdots,A_m)$ such that
\begin{equation}
\mathbb{E}\left[ \|f(x)-f_m(x)\|_{H^2(\Omega)}^2\right] \lesssim \frac{1}{m} \frac{d^{\beta+1}}{C_{\alpha}},	
\end{equation}
where 
\begin{equation}
f_m(x)=\frac{1}{m}\sum\limits_{i=1}^m A_i\sigma(W_i\cdot x+B_i)
\end{equation}
and 
\begin{equation}
\begin{aligned}
A_i= \frac{1}{2\pi \hat{\sigma}(1) }\frac{|\hat{f}(W_i)|\cos(\theta(W_i)-B_i)}{p_1(W_i)p_2(B_i)} I_{\{|B_i|\leq Cd(1+\|W_i\|_2)\} }
\end{aligned}
\end{equation}
for $1\leq i \leq m$.
\end{theorem}

\begin{remark}
Note that $\frac{1}{C_{\alpha}}\sim \frac{1}{\alpha-1} \mathcal{H}^{d-1}(S^{d-1})$, i.e., the hypervolume of the $(d-1)$-dimensional unit sphere. Therefore, when $d$ is sufficiently large, this term becomes small, and $\frac{d^{\beta+1}}{C_{\alpha}}$ does too.
\end{remark}

\begin{remark}
Here, we aim to discuss the relationship between the class of functions we are considering and the Sobolev class of functions. Lemma 2.1 in \cite{41} shows that for function $u$ and any $\epsilon>0$, 
\begin{equation*}
\|u\|_{H^s(\Omega)} \lesssim \|u\|_{\mathcal{B}^s(\Omega)}\lesssim \|u\|_{H^{s+\frac{d}{2}+\epsilon}(\Omega)}.
\end{equation*}
This demonstrates that Sobolev spaces of sufficiently high order can be embedded into Barron spaces. Furthermore, the form of the definition for the functions considered in Condition 1 can be viewed as an equivalent definition of functions in $H^{s/2}(\Omega)$. Therefore, overall, these two cases are not very strict assumptions.

\end{remark}

\subsection{Optimization and Generalization Results}
Within the framework of Physics-Informed Neural Networks (PINNs), we begin with the following elliptic partial differential equation, endowed with Dirichlet boundary conditions.

\begin{equation}	
\left\{
\begin{aligned}
-\Delta u(x)+V(x)u(x)=f(x) & , & in \ \Omega, \\
u(y)=g(y) & , & on \ \partial \Omega,
\end{aligned}
\right.
\end{equation}
where $V,f,g\in L^{\infty}(\Omega)$. 

Let $u^{*}$ be the solution of the given equation (21), we consider two scenarios:

(1) $u^{*}$ belongs to the Barron space $\mathcal{B}^s(\Omega)$;

(2) $u^{*}$ satisfies Condition 1.

In the first scenario, for brevity and readability, we focus on the situation where $s$ is sufficiently large. Specifically, we take $s$ such that $2(s-3)\geq d+1$, i.e., $s\geq\frac{d+7}{2}$. For other values of $s$, as can be seen from the proof in the appendix of \cite{42}, the process does not undergo significant changes; the only adjustment required is in the final step concerning the selection of values for $m$ and $\lambda$.

According to Theorem 1, we can formulate a random neural network $u$ given by the expression:
\begin{equation}
u(x)=\sum\limits_{i=1}^m a_i\sigma(\omega_i \cdot x+b_i),
\end{equation}
where $\omega_1,\cdots,\omega_m\sim P_1$, $b_1,\cdots,b_m\sim P_2$ (as defined in Theorem 1) and $a=(a_1,\cdots,a_m)\in \mathbb{R}^m$ is to be determined. For $1\leq i \leq m$, we denote the $k$-th component of $\omega_i$ by $\omega_{i,k}$.

Given the form of the random neural network $u$, the loss function of PINNs can be written as
\begin{equation}
\begin{aligned}
&L(a):= \int_{\Omega} 	(-\Delta u+Vu-f)^2dx + \int_{\partial \Omega} (u-g)^2dy \\
&= \int_{\Omega}(\sum\limits_{i=1}^m a_i(-\sum\limits_{k}^d  \omega_{i,k}^2\sigma^{''}(\omega_i\cdot x+b_i)+V(x)\sigma(\omega_i \cdot x+b_i)     ) \\
&\quad -f(x))^2dx+ \int_{\partial \Omega}(\sum\limits_{i=1}^m a_i \sigma(\omega_i \cdot y+b_i)-g(y))^2dy \\
&= \int_{\Omega} (a\cdot f_1(x)+g_1(x))^2 dx+ \int_{\partial \Omega} (a\cdot f_2(y)+g_2(y))^2dy,
\end{aligned}
\end{equation}
where $g_1(x)=-f(x), g_2(y)=-g(y)$ and the $i$-th components of $\mathbb{R}^m$-valued functions $f_1(x)$ and $f_2(y)$ are defined as
\begin{equation}
\begin{aligned}
f_1^i(x)&=-\sum\limits_{k=1}^d \omega_{i,k}^2 \sigma^{''}(\omega_i \cdot x+b_i) +V(x)\sigma(\omega_i \cdot x+b_i) \\
&= |\omega_i|^2\sigma^{''}(\omega_i \cdot x+b_i) +V(x)\sigma(\omega_i \cdot x+b_i)
\end{aligned}
\end{equation}
and $f_2^i(y)= \sigma(\omega_i\cdot y+b_i)$, respectively.

To solve for $a$, we introduce a regularization term to the loss function $L(\cdot)$, resulting in the following objective function:
\begin{equation}
\begin{aligned}
&F(a):=L(a)+\lambda\|a\|_2^2 \\
&=\int_{\Omega} (a\cdot f_1(x)+g_1(x))^2 dx+ \int_{\partial \Omega} (a\cdot f_2(y)+g_2(y))^2dy \\
&\quad +\lambda\|a\|_2^2,
\end{aligned}
\end{equation} 
where $\lambda$ is a hyperparameter that controls the regularization strength and is to be determined later.

For simplicity, assuming an equal number of samples from the interior and the boundary, the empirical version of $F(a)$ can be expressed as
\begin{equation}
F_n(a):=\frac{1}{n}\sum\limits_{i=1}^n l(a;x_i,y_i)+\lambda \|a\|_2^2,
\end{equation}
where
\begin{equation}
l(a;x_i,y_i)=|\Omega|(a\cdot f_1(x_i)+g_1(x_i))^2+|\partial \Omega|(a\cdot f_2(y_i)+g_2(y_i))^2.
\end{equation}

In the first case, according to Theorem 1, there exists a random neural network $u_m=\sum\limits_{i=1}^m \bar{a}_i\sigma(\omega_i\cdot x+b_i)$ such that 
\begin{equation}
\mathbb{E}\left[\|u(x)-u_m(x)\|_{H^2(\Omega)}^2 \right] \lesssim \frac{1}{\sqrt{m}},
\end{equation}
where 
\begin{equation}
\bar{a}=(\bar{a}_1,\cdots,\bar{a}_m), \ |\bar{a}_i|\lesssim 1/\sqrt{m}
\end{equation}
for $1\leq i\leq m$.

Note that $\|\bar{a}\|_2=\sqrt{\sum_{i=1}^m\bar{a}_i^2 }=\mathcal{O}(1)$, thus we can employ the projected gradient descent method for 
\begin{equation}
a\in \mathcal{A}:=\{a\in \mathbb{R}^m : \|a\|_2\leq C \},
\end{equation}
where $C$ is a constant ensuring that $\bar{a}\in \mathcal{A}$. Specifically, the method of projected gradient descent consists of iteration of the following update rule for $t=1,\cdots,T,$
\begin{equation}
\left\{	
\begin{aligned}
y_t &= a_t-\nu_t g_t, \quad where \ g_t\in \partial F_n(a_t)\\
a_{t+1}&=P_{\mathcal{A}}(y_t),
\end{aligned}
\right.
\end{equation}
where $T$ is the total number of iterations, $a_1$ is an initial approximation and $P_{\mathcal{A}}$ is the projection operator onto the convex closed set $\mathcal{A}$. The step size $\nu_t$ is a positive scalar chosen according to a specific schedule or rule.

From the specific form of $F_n(\cdot)$, we can see that it is strongly convex and smooth on $\mathcal{A}$. In convex optimization, it is well-known that for a function $f$ that is $\alpha$-strongly convex and $\beta$-smooth on $\mathcal{A}$, the projected gradient descent with $\nu_t=1/\beta$ ensures that for all non-negative integers $t$, the following inequality holds for the sequence $\{y_t\}$ generated by (31):
\begin{equation}
\|y_{t+1}-y^{*}\|_2^2 \leq \exp{\left(-\frac{t}{\kappa}\right)} \|y_1-y^{*}\|_2^2,
\end{equation}
where $\kappa=\beta/\alpha$. 

The strong convexity of $F_n(\cdot)$ is crucial for the optimization process, ensuring that the algorithm converges to the optimal solution efficiently. Moreover, this property contributes to the generalization performance of the random neural networks, which can achieve the fast rate $\mathcal{O}(1/n)$. 

\begin{theorem}
Assume that $u^{*}\in \mathcal{B}^s(\Omega)$ with $s\geq\frac{d+7}{2}$, by applying the projected gradient descent in (31) with $\nu_t=\mathcal{O}(m^{-\frac{d+3}{d+1}})$ to the empirical loss function $F_n(\cdot)$, we have that
\begin{equation}
\begin{aligned}
&\mathbb{E}_{W,B}\mathbb{E}_{X,Y}[L(a_T;f_1,f_2,g_1,g_2)] \\
&\lesssim \frac{1}{\sqrt{m}}+\lambda+\frac{m^{3+\frac{6}{d+1}}}{\lambda^2}e^{-\frac{\lambda T}{m^{1+2/(d+1)}}}+\frac{m^{2+\frac{4}{d+1}}}{\lambda n}\log\log(mn),
\end{aligned}
\end{equation}
where
\begin{equation}
\begin{aligned}
&L(a_T;f_1,f_2,g_1,g_2):=\\
&=\int_{\Omega} (a_T \cdot f_1(x)+g_1(x))^2dx+\int_{\partial\Omega} (a_T \cdot f_2(y)+g_2(y))^2dy
\end{aligned}
\end{equation}
and $a_T$ is obtained from (31) after $T$ iterations.

Taking $\lambda=\frac{1}{\sqrt{m}}$, $m=n^{\frac{d+1}{3d+7}}$ and $T=\mathcal{O}(\sqrt{n}\log(n))$, we have
\begin{equation}
\mathbb{E}_{W,B}\mathbb{E}_{X,Y}[L(a_T;f_1,f_2,g_1,g_2)]=\mathcal{O}(n^{-\frac{d+1}{2(3d+7)} }\log\log(n)).
\end{equation}
\end{theorem}

The first term represents the approximation error. The second term arises from the addition of the regularization term, which imposes strong convexity to the empirical loss function. This modification allows for the application of the projected gradient descent algorithm, ensuring rapid convergence. The third term corresponds to the optimization error, while the fourth term is the generalization error. 
The regularization term also contributes to the generalization error, leading to the fast rate $\mathcal{O}(1/n)$, which is attributed to the strong convexity of the empirical loss function and the utilization of the localization technique.

\begin{remark}
In Theorem 3, we have not show the dependency on the dimension, which is an important factor. In fact, from the proof, the final bound provided in Theorem 1 depends quadratically on $|\Omega|$ and $|\partial \Omega|$. Furthermore, the bound is linear with respect to the constant involved in the Sobolev trace theorem, which is utilized for estimating the approximation error. The generalization error can also be quantified by an appropriate Sobolev norm. Specifically, as demonstrated in \cite{4}, we have 
\begin{equation}
\|u-u^{*}\|_{H^{\frac{1}{2}}(\Omega)}^2 \lesssim \|-\Delta u+Vu-f\|_{L^2(\Omega)}^2+\|u-g\|_{L^2(\partial \Omega)}^2,
\end{equation}
where $u^{*}$ is the solution of equation (21). 
\end{remark}

\begin{remark}
In \cite{2}, algorithm stability has been utilized to derive the excess risk for the projected gradient descent method applied to a loss function that is $\lambda$-strongly convex and $L$-Lipschitz continuous. It has been demonstrated that $T=\mathcal{O}(n^2)$ iterations are required to attain the excess risk $\tilde{\mathcal{O}}(\frac{L^2}{\lambda n})$. It is important to note that our approach to deriving the generalization bound relies solely on the strong convexity and Lipschitz continuity of the loss function. Under the same conditions as in \cite{2}, by employing the projected gradient descent method with $\nu_t=\frac{2}{\lambda(t+1)}$ and $a_T=\frac{2}{T(T+1)}\sum_{t=1}^{T}ta_t$, our method can achieve the same excess risk within $T=\mathcal{O}(n)$ steps. Additionally, because the loss function in our setting is also smooth, $T=\mathcal{O}(\sqrt{n}\log(n))$ is sufficient to reach the same excess risk.
\end{remark}

In the second case, because the randomly sampled hidden weights are no longer bounded, the loss function can become unbounded. Consequently, projected gradient descent is not suitable for this optimization problem, as the output weights of the random neural networks need to be bounded to achieve the approximation results stated in Theorem 2. Moreover, standard gradient descent is also unsuitable. Although the empirical loss function is convex, it lacks Lipschitz continuity, which means there is no guarantee of convergence for gradient descent. Instead, we can directly apply the closed form of ridge regression. According to the following theorem, we can choose $m=n^{1/4}$, which provides an acceptable level of computational complexity for employing closed-form solutions.

\begin{theorem}
Assume that $u^{*}$ satisfies Condition 1, and let $a_n=\mathop{\arg\min}_{a\in \mathbb{R}^m} F_n(a)$, then we can deduce that 
\begin{equation}
\mathbb{E}_{W,B}\mathbb{E}_{X,Y}[L(a_n;f_1,f_2,g_1,g_2)] \lesssim \frac{1}{m}+ \frac{\lambda}{m}+\frac{m}{\lambda\sqrt{n}},
\end{equation}
where 
\begin{equation}
\begin{aligned}
&L(a_n;f_1,f_2,g_1,g_2):=\\
&=\int_{\Omega} (a_n \cdot f_1(x)+g_1(x))^2dx+\int_{\partial\Omega} (a_n \cdot f_2(y)+g_2(y))^2dy.
\end{aligned}
\end{equation}
By taking $\lambda=\mathcal{O}(1), m=n^{1/4}$, we have
\begin{equation}
\mathbb{E}_{W,B}\mathbb{E}_{X,Y}[L(a_n;f_1,f_2,g_1,g_2)]=\mathcal{O}(n^{-\frac{1}{4}}).
\end{equation}
Moreover, $a_n$ can be solved by using the closed form of the ridge regression and the complexity is $\mathcal{O}(m^2n)=\mathcal{O}(n^{\frac{3}{2}})$.
\end{theorem}

\begin{remark}
The error components are distinctly categorized: the first represents the approximation error, the second is attributed to regularization, and the third captures the generalization error. In Theorem 4, an optimization error is notably absent, contrasting with Theorem 3, due to the implementation of the closed-form solution from ridge regression. Furthermore, the unbounded nature of the loss function makes it impossible to to obtain a fast rate for generalization error.
\end{remark}

\section{Experimental Results}
\subsection{Verification of Theorem}
Theorem 4 shows if $m=n^{\frac{1}{4}}$, where $m$ is the hidden width and $n$ is the collocation points number, then the testing loss is expected to be $\mathcal{O}(\frac{1}{m})$.
We demonstrate this result from a simple 1D Poisson equation.

\begin{figure}[htbp]
\centering
\includegraphics[width=0.48\textwidth]{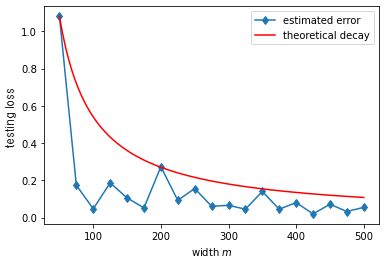}
\caption{Testing loss for the Poisson equation. A random neural network with $m$ hidden nodes ($m$ ranges from 50 to 500) is used to solve 1D Poisson equation with PINN method. The dots show the estimated loss $L(m)$, the line shows the decay implied by the theoretical results $\frac{50L(50)}{m}$.}
\label{loss-poisson}
\end{figure}

The hidden layer's width ranges from 50 to 500, with $tanh$ activations.
In Figure~\ref{loss-poisson}, we can see that the estimated physics-informed loss decays as the width $m$ increases.
The decay trend is consistent with the theoretical decay $\mathcal{O}(\frac{1}{m})$, and the fluctuations may come from the different random initialization of the first layer.

\subsection{Application}
We apply the random neural network to solve the Reaction-Diffusion (RD) equation, which is significant to nonlinear physics, chemistry, and developmental biology. We consider the spatial-temporal Reaction-Diffusion equation with the following form:
\begin{align*}
	\left\{
	\begin{array}{l}
		u_t = d_1 u_{xx} + d_2 u^2, \quad t\in[0,1],x\in[-1,1], \\
		u(0,x) =  \sin(2\pi x)(1+\cos(2\pi x)), \\
		u(t,-1) = u(t, 1)=0,
	\end{array}
	\right.
\end{align*}
where $d_1=d_2=0.01$.

\begin{figure}[htbp]
\centering
\includegraphics[width=0.50\textwidth]{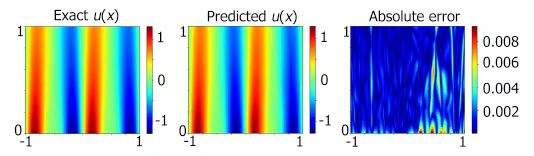}
\caption{Comparison between the reference and predicted solutions for the Reaction-Diffusion equation.}
\label{results-rd}
\end{figure}

We choose $N_{b}=300$ randomly sampled points on the initial domain and the boundary domain, and $N_{f}=5,000$ randomly sampled points within the domain $\Omega=[0,1]\times[-1,1]$. A random neural network with $4$ hidden layers, each containing $128$ units with $\tanh$ activation functions, is used for all computations. The training process runs for $10,000$ epochs with a learning rate of $lr=1e-3$. The relative $L_{2}$ error is $0.13 \%$. As shown in Figure~\ref{results-rd}, the random neural network method achieves high accuracy.
The training time for conventional PINN method is 1397s in \cite{40}, while the training time for random neural network is 3.2s. 
This highlights the high efficiency and precision of our approach.

\section{Conclusion and Discussion}
In this paper, we have established approximation results for Barron-type functions using random neural networks. Our approach diverges from prior studies by providing precise forms for the sampling distributions of hidden weights, rather than relying on empirical selection. Within the PINNs framework, applying random neural networks to solve PDEs translates into formulating regularized least square problems. We address these problems by employing projected gradient descent for one case and leveraging the closed-form solutions for regularized least squares in the other. Combining with the approximation results, we present a comprehensive error analysis. Nevertheless, this work has its limitations. For example, the approximation rate presented in Theorem 1 for functions in $\mathcal{B}^s(\Omega)$ suffers the curse of dimensionality. Moreover, both approximation results in Theorem 1 and Theorem 2 are only valid for relatively small function spaces. Expanding the approximation capabilities of random neural networks to more general functions presents an interesting direction for future research. Furthermore, the methods we employed for deriving the optimization and generalization errors could be extended to other methods involving random neural networks for solving PDEs, such as the Deep Ritz Method.

\section*{Acknowledgments}
This work was partially supported by the National Natural Science Foundation of China (No. 12025104, No. 62106103, ), and the basic research project (ILF240021A24). This work is also partially supported by High Performance Computing Platform of Nanjing University of Aeronautics and Astronautics.

\bibliographystyle{named}
\bibliography{random}


\newpage
\appendix
\onecolumn

\section*{Appendix}

\section{Proof of Theorem 1}

\begin{proof}
Note that the activation function $\sigma$ satisfies the polynomial decay condition, i.e., 
\begin{equation}
\begin{aligned}
|\sigma^{(k)}(t)|&\leq C_p(1+|t|)^{-p}
\end{aligned}
\end{equation}
for $0\leq k \leq 2$ and any $p>1$, as $\sigma$ is compactly supported in $[-2,2]$. In the following, we only use this condition with $p=2$.

Thus for any $x\in \Omega$ and $s\in \mathbb{N}$, we have
\begin{equation}
\begin{aligned}
&\int_{\mathbb{R}} |\sigma^{(k)}(\omega\cdot x+b)|^sdb \\
& \leq C_2\int_{\mathbb{R}}(1+|\omega\cdot x+b|)^{-2s}db \\
&\leq C_2\int_{\mathbb{R}}(1+\max(0, |b|-\|\omega\|_1))^{-2s}db\\
&\leq C_2\int_{\{|b|\leq \|\omega\|_1\}} db + 2C_2\int_{\{b> \|\omega\|_1\}}(1+b-\|\omega\|_1))^{-2s}db \\
&=C_2(2\|\omega\|_1 +2 \frac{1}{2s-1})\\
&\leq 2C_2(1+\|\omega\|_1).
\end{aligned}
\end{equation}

In the integral representation of $f$, we simultaneously truncate both $\omega$ and $b$.
\begin{equation}
\begin{aligned}
f(x)&=\int_{\mathbb{R}^d} \int_{\mathbb{R}} \frac{1}{2\pi \hat{\sigma}(1) }\sigma(\omega\cdot x+b)|\hat{f}(\omega)|\cos(\theta(\omega)-b)dbd\omega \\
&=\int_{\{\|\omega\|_1\leq M,|b|\leq 2M\}}  \frac{1}{2\pi \hat{\sigma}(1) }\sigma(\omega\cdot x+b)|\hat{f}(\omega)|\cos(\theta(\omega)-b)dbd\omega \\
&\quad + \int_{\{ \|\omega\|_1\leq M, |b|\leq 2M\}^c  } \frac{1}{2\pi \hat{\sigma}(1) }\sigma(\omega\cdot x+b)|\hat{f}(\omega)|\cos(\theta(\omega)-b)dbd\omega\\
&:= f_1(x)+f_2(x), 
\end{aligned}
\end{equation}
where we write $\hat{f}$ for $\hat{f_e}$ for simplicity.

From the forms of $f_1(x)$ and $f_2(x)$, we will use a random neural network to approximate $f_1(x)$ via the Monte Carlo method. To ensure that $f(x)$ can also be well-approximated by this random neural network, we need to make sure that $f_2(x)$ has a sufficiently small $H^2$ norm when $M$ is large enough.

For $f_2(x)$, it can be decomposed into two terms: one corresponding to integral over $\{ \|\omega\|_1>M, b\in \mathbb{R}\}$ and anther to the integral over $\{\|\omega\|_1\leq M, |b|>2M \}$.
\begin{equation}
\begin{aligned}
f_2(x)&= \int_{\{ \|\omega\|_1\leq M, |b|\leq 2M\}^c  } \frac{1}{2\pi \hat{\sigma}(1) }\sigma(\omega\cdot x+b)|\hat{f}(\omega)|\cos(\theta(\omega)-b)dbd\omega\\
&= \int_{\{\|\omega\|_1>M \}} \int_{\mathbb{R}}\frac{1}{2\pi \hat{\sigma}(1) }\sigma(\omega\cdot x+b)|\hat{f}(\omega)|\cos(\theta(\omega)-b)dbd\omega \\  &\quad +\int_{\{\|\omega\|_1\leq M,|b|>2M\}} \frac{1}{2\pi \hat{\sigma}(1) }\sigma(\omega\cdot x+b)|\hat{f}(\omega)|\cos(\theta(\omega)-b)dbd\omega\\
&:= f_{21}(x)+f_{22}(x).
\end{aligned}
\end{equation}

In the following, we first estimate the $H^2$ norm of $f_2(x)$. Note that when $\|\omega\|_1\leq M, |b|>2M$, we have 
\[|\omega\cdot x+b|\geq |b|-|\omega\cdot x|\geq |b|-\|\omega\|_1> M \geq 2,\]
which implies that $\sigma(\omega\cdot x+b)=0$, since $\sigma$ is compactly supported in $[-2,2]$. Thus, $f_{22}(x)=0$.

For the estimation of the $H^2$ norm of $f_2(x)$, it remains only to estimate the $H^2$ norm of $f_{21}(x)$, which can be expressed explicitly as follows.
\begin{equation}
\begin{aligned}
\|f_{21}(x)\|_{H^2(\Omega)}^2 &= \int_{\Omega} \left(\int_{\{\|\omega\|_1>M \}} \int_{\mathbb{R}} \frac{1}{2\pi \hat{\sigma}(1) }\sigma(\omega\cdot x+b)|\hat{f}(\omega)|\cos(\theta(\omega)-b)dbd\omega\right)^2dx \\
& \quad + \int_{\Omega}\sum\limits_{k=1}^d \left(\int_{\{\|\omega\|_1>M \}} \int_{\mathbb{R}} \frac{1}{2\pi \hat{\sigma}(1) }\omega_k\sigma^{'}(\omega\cdot x+b)|\hat{f}(\omega)|\cos(\theta(\omega)-b)dbd\omega\right)^2dx \\
&\quad +\int_{\Omega} \sum\limits_{i=1}^d \sum\limits_{j=1}^d  \left(\int_{\{\|\omega\|_1>M \}} \int_{\mathbb{R}} \frac{1}{2\pi \hat{\sigma}(1) }\omega_i \omega_j\sigma^{''}(\omega\cdot x+b)|\hat{f}(\omega)|\cos(\theta(\omega)-b)dbd\omega\right)^2dx\\
&:= \bar{I}_1+\bar{I}_2+\bar{I}_3.
\end{aligned}
\end{equation}

From (41), for $\bar{I}_1$, we have that for any $x\in \Omega$, 
\begin{equation}
\begin{aligned}
&\left(\int_{\{\|\omega\|_1>M \}} \int_{\mathbb{R}}\sigma(\omega\cdot x+b)|\hat{f}(\omega)|\cos(\theta(\omega)-b)dbd\omega\right)^2\\
&\leq \left(\int_{\{\|\omega\|_1>M \}} \int_{\mathbb{R}}|\sigma(\omega\cdot x+b)||\hat{f}(\omega)|dbd\omega\right)^2\\
&\lesssim \left(\int_{\{\|\omega\|_1>M \}} (1+\|\omega\|_1)|\hat{f}(\omega)|d\omega\right)^2\\
&\leq \left(\frac{1}{M^{s-1}}\|f\|_{\mathcal{B}^s(\Omega)}\right)^2,
\end{aligned}
\end{equation}
where the second inequality follows from the polynomial decay condition of $\sigma$.

As for $\bar{I}_2$, for any $k\in \{1,\cdots,d\}$, we can deduce from (41) that
\begin{equation}
\begin{aligned}
&\left(\int_{\{\|\omega\|_1>M \}} \int_{\mathbb{R}}\omega_k\sigma^{'}(\omega\cdot x+b)|\hat{f}(\omega)|\cos(\theta(\omega)-b)dbd\omega\right)^2\\
&\leq \left(\int_{\{\|\omega\|_1>M \}} \int_{\mathbb{R}}|\omega_k||\sigma^{'}(\omega\cdot x+b)||\hat{f}(\omega)|dbd\omega\right)^2\\
&\lesssim \left(\int_{\{\|\omega\|_1>M \}} |\omega_k|(1+\|\omega\|_1)|\hat{f}(\omega)|d\omega\right)^2\\
&\leq \left(\int_{\{\|\omega\|_1>M \}} |\omega_k|^2|\hat{f}(\omega)|d\omega\right) \left(\int_{\{\|\omega\|_1>M \}} (1+\|\omega\|_1)^2|\hat{f}(\omega)|d\omega\right),
\end{aligned}
\end{equation}
where the last inequality follows from the Cauchy inequality.

Summing $k$ from $1$ to $d$ in (46) yields that 
\begin{equation}
\begin{aligned}
\bar{I}_2&=\int_{\Omega}\sum\limits_{k=1}^d \left(\int_{\{\|\omega\|_1>M \}} \int_{\mathbb{R}} \frac{1}{2\pi \hat{\sigma}(1) }\omega_k\sigma^{'}(\omega\cdot x+b)|\hat{f}(\omega)|\cos(\theta(\omega)-b)dbd\omega\right)^2dx \\
&\lesssim  \left(\int_{\{\|\omega\|_1>M \}} (1+\|\omega\|_1)^2|\hat{f}(\omega)|d\omega\right) \left(\int_{\{\|\omega\|_1>M \}} (1+\|\omega\|_1)^2|\hat{f}(\omega)|d\omega\right)\\
&\leq \left( \frac{1}{M^{s-2}} \int_{\{\|\omega\|_1>M \}} (1+\|\omega\|_1)^s|\hat{f}(\omega)|d\omega \right)^2 \\
&\leq \left(\frac{1}{M^{s-2}} \|f\|_{\mathcal{B}^s(\Omega)} \right)^2.
\end{aligned}
\end{equation}

Similarly, for $\bar{I}_3$, we can deduce that for any $i,j$ ($1\leq i, j\leq d$, 
\begin{equation}
\begin{aligned}
&\left(\int_{\{\|\omega\|_1>M \}} \int_{\mathbb{R}} \omega_i \omega_j\sigma^{''}(\omega\cdot x+b)|\hat{f}(\omega)|\cos(\theta(\omega)-b)dbd\omega\right)^2\\
&\lesssim \left(\int_{\{\|\omega\|_1>M \}} |\omega_i \omega_j| (1+\|\omega\|_1)|\hat{f}(\omega)|d\omega\right)^2\\
&\leq  \left(\int_{\{\|\omega\|_1>M \}} |\omega_i|^2  (1+\|\omega\|_1)|\hat{f}(\omega)|d\omega\right) \left(\int_{\{\|\omega\|_1>M \}} |\omega_j|^2  (1+\|\omega\|_1)|\hat{f}(\omega)|d\omega\right),
\end{aligned}
\end{equation}
where the first inequality follows from (41) and the second inequality is from the Cauchy inequality.

Therefore, 
\begin{equation}
\begin{aligned}
\bar{I}_3&=\sum\limits_{i=1}^d\sum\limits_{j=1}^d\left(\int_{\{\|\omega\|_1>M \}} \int_{\mathbb{R}} \frac{1}{2\pi \hat{\sigma}(1) }\omega_i \omega_j\sigma^{''}(\omega\cdot x+b)|\hat{f}(\omega)|\cos(\theta(\omega)-b)dbd\omega\right)^2\\
&\lesssim \left(\sum\limits_{i=1}^d\int_{\{\|\omega\|_1>M \}} |\omega_i|^2  (1+\|\omega\|_1)|\hat{f}(\omega)|d\omega\right) \left(\sum\limits_{j=1}^d\int_{\{\|\omega\|_1>M \}} |\omega_j|^2  (1+\|\omega\|_1)|\hat{f}(\omega)|d\omega\right) \\
&\leq \left(\int_{\{\|\omega\|_1>M \}}  (1+\|\omega\|_1)^3|\hat{f}(\omega)|d\omega\right)^2 \\
&\leq \left(  \frac{1}{M^{s-3}} \|f\|_{\mathcal{B}^s(\Omega)} \right)^2.
\end{aligned}
\end{equation}

Combining these estimations in (45), (47) and (49), we have 
\begin{equation}
\begin{aligned}
\|f_{21}(x)\|_{H^2(\Omega)}^2 &= \bar{I}_1+\bar{I}_2+\bar{I}_3 \lesssim \frac{\|f\|_{\mathcal{B}^s(\Omega)}}{M^{s-3}}.
\end{aligned}
\end{equation}

It remains to consider the stochastic approximation of the main term $f_1(x)$, recall that 
\begin{equation}
f_1(x)=\int_{\{\|\omega\|_1\leq M,|b|\leq 2M\}}  \frac{1}{2\pi \hat{\sigma}(1) }\sigma(\omega\cdot x+b)|\hat{f}(\omega)|\cos(\theta(\omega)-b)dbd\omega.
\end{equation}

Thus, we can express $f_1(x)$ in the form of expectation.
\begin{equation}
\begin{aligned}
f_1(x)&=\int_{\{\|\omega\|_1\leq M,|b|\leq 2M\}}  \frac{1}{2\pi \hat{\sigma}(1) }\frac{\sigma(\omega\cdot x+b)|\hat{f}(\omega)|\cos(\theta(\omega)-b)}{p_1(\omega)p_2(b)}  p_1(\omega) p_2(b)dbd\omega\\
&= \mathbb{E}_{P_1(\omega), P_2(b)}\left[ \frac{1}{2\pi \hat{\sigma}(1) }\frac{\sigma(\omega\cdot x+b)|\hat{f}(\omega)|\cos(\theta(\omega)-b)}{p_1(\omega)p_2(b)}\right] ,
\end{aligned}
\end{equation}
where $P_1(\omega):=p_1(\omega)d\omega, P_2(b):=p_2(b)db$ are uniform distributions over $\{\|\omega\|_1\leq M\}$ and $\{|b|\leq 2M\}$ respectively, i.e., $p_1(\omega)=\frac{1}{C_dM^d}, p_2(b)=\frac{1}{4M}$, where $C_d$ is the measure of $d$-dimensional unit $L^1$ ball.

Thus, $f_1(x)$ can also be written as
\begin{equation}
f_1(x)= \mathbb{E}_{P_1(W), P_2(B)} \left[ A\sigma(W\cdot x+B) \right],\
\end{equation}
where 
\begin{equation}
A= \frac{1}{2\pi \hat{\sigma}(1) }\frac{|\hat{f}(W)|\cos(\theta(W)-B)}{p_1(W)p_2(B)} .
\end{equation}

From the form of $f_1(x)$ in (53), we can take $W_1,\cdots, W_m \sim P_1$ and $B_1,\cdots, B_m\sim P_2$, then let 
\begin{equation}
A_i= \frac{1}{2\pi \hat{\sigma}(1) }\frac{|\hat{f}(W_i)|\cos(\theta(W_i)-B_i)}{p_1(W_i)p_2(B_i)} , 1\leq i \leq m.
\end{equation}

Thus, we can construct the random neural network $f_m(x)$:
\begin{equation}
f_m(x)= \frac{1}{m}\sum\limits_{i=1}^m A_i\sigma(W_i\cdot x+B_i).
\end{equation}

Finally, it is sufficient to estimate $\mathbb{E}\left[\|f_1(x)-f_m(x)\|_{H^2(\Omega)}^2\right]$, which can be expressed as follows.
\begin{equation}
\begin{aligned}
&\mathbb{E}\left[\|f_1(x)-f_m(x)\|_{H^2(\Omega)}^2\right]\\
&= \mathbb{E}\left[  \left\|f_1(x)- \frac{1}{m} \sum\limits_{i=1}^m A_i\sigma(W_i\cdot x+B_i)\right\|_{H^2(\Omega)}^2 \right] \\
&= \mathbb{E}\left[\int_{\Omega}  \left( f_1(x)- \frac{1}{m} \sum\limits_{i=1}^m A_i\sigma(W_i\cdot x+B_i)\right)^2dx  \right]\\
&\quad + \mathbb{E}\left[\sum\limits_{k=1}^d\int_{\Omega}   \left( \partial_{k}f_1(x)- \frac{1}{m} \sum\limits_{i=1}^m A_iW_{i,k}\sigma^{'}(W_i\cdot x+B_i)\right)^2dx    \right] \\
&\quad + \mathbb{E}\left[\sum\limits_{s=1}^d \sum\limits_{t=1}^d\int_{\Omega}   \left( \partial_{s,t}f_1(x)- \frac{1}{m} \sum\limits_{i=1}^m A_iW_{i,s}W_{i,t}\sigma^{''}(W_i\cdot x+B_i)\right)^2dx    \right]\\
&:=I_1+I_2+I_3.
\end{aligned}
\end{equation}

For the first term $I_1$, applying Fubini's theorem yields that
\begin{equation}
\begin{aligned}
&\mathbb{E}\left[\int_{\Omega}  \left( f_1(x)- \frac{1}{m} \sum\limits_{i=1}^m A_i\sigma(W_i\cdot x+B_i)\right)^2dx  \right] \\
&= 	\int_{\Omega} \mathbb{E}\left[  \left( f_1(x)- \frac{1}{m} \sum\limits_{i=1}^m A_i\sigma(W_i\cdot x+B_i)\right)^2  \right]dx \\
&= \int_{\Omega} \frac{1}{m}Var(A\sigma(W\cdot x+B) ) dx\\
&\leq \frac{1}{m} \int_{\Omega} \mathbb{E}[(A\sigma(W\cdot x+B))^2] dx \\
&=\frac{1}{m} \int_{\Omega}  \int_{ \{\|\omega\|_1\leq M\}}  \int_{ \{|b|\leq 2M\}}  \left(\frac{1}{2\pi \hat{\sigma}(1) }\frac{\sigma(\omega\cdot x+b)|\hat{f}(\omega)|\cos(\theta(\omega)-b)}{p_1(\omega)p_2(b)}\right)^2 p_1(\omega) p_2(b)dbd\omega        dx\\
&\lesssim \frac{1}{m} C_dM^{d+1} \int_{\{\|\omega\|_1\leq M\}}  (1+\|\omega\|_1) |\hat{f}(\omega)|^2 d\omega \\
&\lesssim \frac{1}{m} C_dM^{d+1} \|\hat{f}\|_{L^{\infty}(\mathbb{R}^d)}\|f\|_{\mathcal{B}^1(\Omega)},
\end{aligned}
\end{equation}
where the second inequality follows from (41).

Similarly, for the second term $I_2$, we have that for any $k\in \{1,\cdots, d\}$,
\begin{equation}
\begin{aligned}
&\mathbb{E}\left[\int_{\Omega}   \left( \partial_{k}f_1(x)- \frac{1}{m} \sum\limits_{i=1}^m A_iW_{i,k}\sigma^{'}(W_i\cdot x+B_i)\right)^2dx    \right] \\
&\leq \frac{1}{m} \int_{\Omega} Var(AW_k\sigma^{'}(W\cdot x+B) )dx \\
&\leq \frac{1}{m} \int_{\Omega} \mathbb{E}[(AW_k\sigma^{'}(W\cdot x+B) )^2]dx\\
&=  \frac{1}{m} \int_{\Omega} \int_{ \{\|\omega\|_1\leq M\}}  \int_{ \{|b|\leq 2M\}}   \left(\frac{1}{2\pi \hat{\sigma}(1) }\frac{\omega_k\sigma^{'}(\omega\cdot x+b)|\hat{f}(\omega)|\cos(\theta(\omega)-b)}{p_1(\omega)p_2(b)}\right)^2 p_1(\omega) p_2(b)dbd\omega        dx\\
&\lesssim \frac{1}{m} \int_{\Omega} \int_{ \{\|\omega\|_1\leq M\}}  \int_{ \{|b|\leq 2M\}}  C_dM^{d+1} \omega_k^2|\sigma^{'}(\omega\cdot x+b)|^2 |\hat{f}(\omega)|^2dbd\omega dx\\
&\lesssim \frac{1}{m} C_dM^{d+1} \int_{\{\|\omega\|_1\leq M\}}  \omega_k^2 (1+\|\omega\|_1)|\hat{f}(\omega)|^2d\omega.
\end{aligned}
\end{equation}
Summing over $k$ yields that
\begin{equation}
\begin{aligned}
& \mathbb{E}\left[\sum\limits_{k=1}^d\int_{\Omega}   \left( \partial_{k}f_1(x)- \frac{1}{m} \sum\limits_{i=1}^m A_iW_{i,k}\sigma^{'}(W_i\cdot x+B_i)\right)^2dx    \right] \\
&\lesssim 	 \frac{1}{m} C_dM^{d+1} \int_{\{\|\omega\|_1\leq M\}}  (1+\|\omega\|_1)^2 (1+\|\omega\|_1)|\hat{f}(\omega)|^2d\omega \\
&\lesssim \frac{1}{m} C_dM^{d+1} \|\hat{f}\|_{L^{\infty}(\mathbb{R}^d)} \|f\|_{\mathcal{B}^3(\Omega)}.
\end{aligned}
\end{equation}
As for the last term $I_3$, we have that for any $i,j$ ($1\leq i,j\leq d$),
\begin{equation}
\begin{aligned}
&\mathbb{E}\left[\int_{\Omega}   \left( \partial_{s,t}f_1(x)- \frac{1}{m} \sum\limits_{i=1}^m A_iW_{i,s}W_{i,t}\sigma^{''}(W_i\cdot x+B_i)\right)^2dx    \right]	 \\
&\lesssim \frac{1}{m}  \int_{\Omega} \mathbb{E}[(AW_{s}W_{t}\sigma^{''}(W\cdot x+B))^2]dx\\
&\lesssim \frac{1}{m} C_dM^{d+1} \int_{\Omega} \int_{ \{\|\omega\|_1\leq M\}}  \int_{ \{|b|\leq 2M\}} \omega_s^2 \omega_t^2 |\sigma^{''}(\omega\cdot x+b)|^2 |\hat{f}(\omega)|^2 dbd\omega dx\\
&\lesssim \frac{1}{m} C_dM^{d+1} \int_{ \{\|\omega\|_1\leq M\}}  \omega_s^2 \omega_t^2(1+\|\omega\|_1)^2|\hat{f}(\omega)|^2d\omega.
\end{aligned}
\end{equation}
Summing over $i,j$ yields that
\begin{equation}
\begin{aligned}
&\mathbb{E}\left[\sum\limits_{s=1}^d \sum\limits_{t=1}^d\int_{\Omega}   \left( \partial_{s,t}f_1(x)- \frac{1}{m} \sum\limits_{i=1}^m A_iW_{i,s}W_{i,t}\sigma^{''}(W_i\cdot x+B_i)\right)^2dx    \right] \\
&\lesssim \frac{1}{m}C_dM^{d+1}  \int_{\{\|\omega\|_1\leq M\}} \|\omega\|_1^4 (1+\|\omega\|_1)^2|\hat{f}(\omega)|^2 d\omega.
\end{aligned}
\end{equation}
Thus, for $3<s\leq 6$, 
\begin{equation}
\int_{\{\|\omega\|_1\leq M\}} \|\omega\|_1^4 (1+\|\omega\|_1)^2|\hat{f}(\omega)|^2 d\omega \leq M^{6-s}\int_{\{\|\omega\|_1\leq M\}}  (1+\|\omega\|_1)^s|\hat{f}(\omega)|^2 d\omega\leq M^{6-s} \|\hat{f}\|_{L^{\infty}(\mathbb{R}^d)}\|f\|_{\mathcal{B}^s(\Omega)}
\end{equation}
and for $s>6$, 
\begin{equation}
\int_{\{\|\omega\|_1\leq M\}} \|\omega\|_1^4 (1+\|\omega\|_1)^2|\hat{f}(\omega)|^2 d\omega \leq \int_{\{\|\omega\|_1\leq M\}}  (1+\|\omega\|_1)^6|\hat{f}(\omega)|^2 d\omega\leq \|\hat{f}\|_{L^{\infty}(\mathbb{R}^d)}\|f\|_{\mathcal{B}^s(\Omega)}.
\end{equation}

Combining the estimations for $f_1(x)$ (i.e., (58), (60), (62), (63) and (64)) and the estimations for $f_2(x)$ (i.e., (50)), we obtain that
\begin{equation}
\mathbb{E}\left[ \|f(x)-f_m(x)\|_{H^2(\Omega)}^2\right] \lesssim \left\{
\begin{aligned}
\frac{M^{d+7-s}}{m} \|\hat{f}\|_{L^{\infty}(\mathbb{R}^d)}\|f\|_{\mathcal{B}^s(\Omega)}+\left(\frac{1}{M^{s-3}} \|f\|_{\mathcal{B}^s(\Omega)} \right)^2 & , & 3<s<6, \\
\frac{M^{d+1}}{m} \|\hat{f}\|_{L^{\infty}(\mathbb{R}^d)}\|f\|_{\mathcal{B}^s(\Omega)}+\left(\frac{1}{M^{s-3}} \|f\|_{\mathcal{B}^s(\Omega)} \right)^2 & , & s\geq 6.
\end{aligned}
\right.
\end{equation}
Moreover, since
\begin{equation}
A_i= \frac{1}{2\pi \hat{\sigma}(1) }\frac{|\hat{f}(W_i)|\cos(\theta(W_i)-B_i)}{p_1(W_i)p_2(B_i)} , 1\leq i \leq m,
\end{equation}
we can deduce that for $1\leq i \leq m$,
\begin{equation}
|A_i|\lesssim \|\hat{f}\|_{L^{\infty}(\mathbb{R}^d)}M^{d+1}.
\end{equation}
\end{proof}

\section{Proof of Theorem 2}
\begin{proof}
Recall that
\begin{equation}
\begin{aligned}
f(x)	&=  \int_{\mathbb{R}^d}\int_{\mathbb{R}} \frac{1}{2\pi \hat{\sigma}(1) }\sigma(\omega\cdot x+b)|\hat{f}(\omega)|\cos(\theta(\omega)-b)dbd\omega\\
&= \int_{\mathbb{R}^d}\int_{\mathbb{R}} \frac{1}{2\pi \hat{\sigma}(1) } \sigma(\omega\cdot x+b)|\hat{f}(\omega)|\cos(\theta(\omega)-b)  \frac{(1+\|\omega\|_2)^{\alpha}(1+|b|)^{\beta}}{C_{\alpha}C_{\beta}} \frac{C_{\alpha}C_{\beta}}{(1+\|\omega\|_2)^{\alpha} (1+|b|)^{\beta}}dbd\omega\\
&=\mathbb{E}_{P_1(\omega), P_2(b)}\left[ \frac{1}{2\pi \hat{\sigma}(1) }\frac{\sigma(\omega\cdot x+b)|\hat{f}(\omega)|\cos(\theta(\omega)-b)}{p_1(\omega)p_2(b)}\right] ,
\end{aligned}
\end{equation}
where $P_1(\omega)=p_1(\omega)d\omega, P_2(b)=p_2(b)db $ and $p_1(\omega):=\frac{C_{\alpha}}{(1+\|\omega\|_2)^{\alpha}}, p_2(b):=\frac{C_{\beta}}{(1+|b|)^{\beta}}$ are density functions in $\mathbb{R}^d$ and $\mathbb{R}$ respectively. 

Taking $W\sim P_1, B\sim P_2$ and let
\begin{equation}
A= \frac{1}{2\pi \hat{\sigma}(1) }\frac{|\hat{f}(W)|\cos(\theta(W)-B)}{p_1(W)p_2(B)} I_{\{|B|\leq Cd(1+\|W\|_2)\} },
\end{equation}
then we can deduce that
\begin{equation}
f(x)= \mathbb{E}_{P_1(W), P_2(B)} \left[ A\sigma(W\cdot x+B) \right],
\end{equation}
since $\sigma$ is compactly supported in $[-2,2]$, which implies 
\begin{equation}
\sigma(W\cdot x+B)=\sigma(W\cdot x+B)I\{|W\cdot x+B|\leq 2\}=\sigma(W\cdot x+B)I\{|B|\leq C(1+\|W\|_1)\}=\sigma(W\cdot x+B)I\{|B|\leq Cd(1+\|W\|_2)\}.
\end{equation}

Therefore, we can construct a random neural network $f_m$ as follows.
\begin{equation}
f_m(x)=\frac{1}{m} \sum\limits_{i=1}^m A_i\sigma(W_i\cdot x+B_i),
\end{equation}
where $W_1,\cdots, W_m \sim P_1$ and $B_1,\cdots, B_m\sim P_2$,
\begin{equation}
A_i= \frac{1}{2\pi \hat{\sigma}(1) }\frac{|\hat{f}(W_i)|\cos(\theta(W_i)-B_i)}{p_1(W_i)p_2(B_i)} I_{\{|B_i|\leq Cd(1+\|W_i\|_2)\} }, 1\leq i \leq m.
\end{equation}

It remains to estimate $\mathbb{E}\left[\|f(x)-f_m(x)\|_{H^2(\Omega)}^2\right]$, which can be expressed explicitly as follows.
\begin{equation}
\begin{aligned}
&\mathbb{E}\left[\|f(x)-f_m(x)\|_{H^2(\Omega)}^2\right]\\
&= \mathbb{E}[  \left\|f(x)- \frac{1}{m} \sum\limits_{i=1}^m A_i\sigma(W_i\cdot x+B_i)\right\|_{H^2(\Omega)}^2  ] \\
&= \mathbb{E}\left[\int_{\Omega}  \left( f(x)- \frac{1}{m} \sum\limits_{i=1}^m A_i\sigma(W_i\cdot x+B_i)\right)^2dx  \right]\\
&\quad + \mathbb{E}\left[\sum\limits_{k=1}^d\int_{\Omega}   \left( \partial_{k}f(x)- \frac{1}{m} \sum\limits_{i=1}^m A_iW_{i,k}\sigma^{'}(W_i\cdot x+B_i)\right)^2dx    \right] \\
&\quad + \mathbb{E}\left[\sum\limits_{s=1}^d \sum\limits_{t=1}^d\int_{\Omega}   \left( \partial_{s,t}f(x)- \frac{1}{m} \sum\limits_{i=1}^m A_iW_{i,s}W_{i,t}\sigma^{''}(W_i\cdot x+B_i)\right)^2dx    \right]\\
&:=\tilde{I}_1+\tilde{I}_2+\tilde{I}_3.
\end{aligned}
\end{equation}

For the first term $\tilde{I}_1$, applying Fubini's theorem yields that
\begin{equation}
\begin{aligned}
I_1&=\mathbb{E}\left[\int_{\Omega}  \left( f(x)- \frac{1}{m} \sum\limits_{i=1}^m A_i\sigma(W_i\cdot x+B_i)\right)^2dx  \right] \\
&= 	\int_{\Omega} \mathbb{E}\left[  \left( f(x)- \frac{1}{m} \sum\limits_{i=1}^m A_i\sigma(W_i\cdot x+B_i)\right)^2  \right]dx \\
&= \int_{\Omega} \frac{1}{m}Var(A\sigma(W\cdot x+B) ) dx\\
&\leq \frac{1}{m} \int_{\Omega} \mathbb{E}[(A\sigma(W\cdot x+B))^2] dx \\
&=\frac{1}{m} \int_{\Omega}  \int_{\mathbb{R}^d}  \int_{\mathbb{R}}  \left(\frac{1}{2\pi \hat{\sigma}(1) }\frac{\sigma(\omega\cdot x+b)|\hat{f}(\omega)|\cos(\theta(\omega)-b)}{p_1(\omega)p_2(b)}\right)^2 p_1(\omega) p_2(b)dbd\omega        dx\\
&\lesssim \frac{1}{m} \frac{1}{C_{\alpha}C_{\beta}} \int_{\Omega}\int_{\mathbb{R}^d} \int_{\mathbb{R}} (1+\|\omega\|_2)^{\alpha}(1+|b|)^{\beta}|\sigma(\omega\cdot x+b)|^2 |\hat{f}(\omega)|^2 dbd\omega dx\\
&\lesssim \frac{1}{m} \frac{1}{C_{\alpha}C_{\beta}}d^{\beta} \int_{\Omega}\int_{\mathbb{R}^d} \int_{\mathbb{R}} (1+\|\omega\|_2)^{\alpha+\beta}|\sigma(\omega\cdot x+b)|^2 |\hat{f}(\omega)|^2 dbd\omega dx \\
&\lesssim \frac{1}{m} \frac{1}{C_{\alpha}C_{\beta}} \int_{\mathbb{R}^d} d^{\beta+1} (1+\|\omega\|_2)^{\alpha+\beta+1} |\hat{f}(\omega)|^2 d\omega,
\end{aligned}
\end{equation}
where the third inequality follows from that $|b|\leq Cd(1+\|\omega\|_2)$ and the last inequality follows from the polynomial decay condition of $\sigma$, i.e., (41).

Similarly, for the second term $I_2$, we have
\begin{equation}
\begin{aligned}
I_2&= \mathbb{E}\left[\sum\limits_{k=1}^d\int_{\Omega}   \left( \partial_{k}f_1(x)- \frac{1}{m} \sum\limits_{i=1}^m A_iW_{i,k}\sigma^{'}(W_i\cdot x+B_i)\right)^2dx    \right] \\
&\leq \frac{1}{m}\sum\limits_{k=1}^d \int_{\Omega} \int_{\mathbb{R}^d}  \int_{\mathbb{R}}  \left(\frac{1}{2\pi \hat{\sigma}(1) }\frac{\omega_k\sigma^{'}(\omega\cdot x+b)|\hat{f}(\omega)|\cos(\theta(\omega)-b)}{p_1(\omega)p_2(b)}\right)^2 p_1(\omega) p_2(b)dbd\omega        dx\\
&\lesssim \frac{1}{m} \frac{1}{C_{\alpha}C_{\beta}} \sum\limits_{k=1}^d \int_{\Omega}\int_{\mathbb{R}^d} \int_{\mathbb{R}} \omega_k^2(1+\|\omega\|_2)^{\alpha}(1+|b|)^{\beta}|\sigma(\omega\cdot x+b)|^2 |\hat{f}(\omega)|^2 dbd\omega dx\\
&\lesssim \frac{1}{m}\frac{1}{C_{\alpha}C_{\beta}}d^{\beta+1} \sum\limits_{k=1}^d\int_{\mathbb{R}^d}  \omega_k^2 (1+\|\omega\|_2)^{\alpha+\beta+1}|\hat{f}(\omega)|^2d\omega\\
&= \frac{1}{m}\frac{1}{C_{\alpha}C_{\beta}}d^{\beta+1} \int_{\mathbb{R}^d}  \|\omega\|_2^2 (1+\|\omega\|_2)^{\alpha+\beta+1}|\hat{f}(\omega)|^2d\omega\\
&\lesssim \frac{1}{m}\frac{1}{C_{\alpha}C_{\beta}}d^{\beta+1} \int_{\mathbb{R}^d}  (1+\|\omega\|_2)^{\alpha+\beta+3}|\hat{f}(\omega)|^2d\omega.
\end{aligned}
\end{equation}

As for the third term $I_3$, we can deduce that
\begin{equation}
\begin{aligned}
I_3&=\mathbb{E}\left[\sum\limits_{s=1}^d \sum\limits_{t=1}^d\int_{\Omega}   \left( \partial_{s,t}f_1(x)- \frac{1}{m} \sum\limits_{i=1}^m A_iW_{i,s}W_{i,t}\sigma^{''}(W_i\cdot x+B_i)\right)^2dx    \right] \\
&\lesssim \frac{1}{m}\frac{1}{C_{\alpha}C_{\beta}}d^{\beta+1}\int_{\mathbb{R}^d} (1+\|\omega\|_2)^{\alpha+\beta+5}|\hat{f}(\omega)|^2d\omega,
\end{aligned}
\end{equation}
which is the reason of the requirement that $s>\alpha+\beta+5$.

Finally, combining (75), (76) and (77) leads to the conclusion:
\begin{equation}
\mathbb{E}\left[ \|f(x)-f_m(x)\|_{H^2(\Omega)}^2\right] \lesssim \frac{1}{m} \frac{d^{\beta+1}}{C_{\alpha}}.
\end{equation}
\end{proof}

\section{Proof of Theorem 3}
\begin{proof}
Recall that the loss function and its empirical part are 
\begin{equation}
F(a)=L(a)+\lambda\|a\|_2^2=\int_{\Omega} (a\cdot f_1(x)+g_1(x))^2 dx+ \int_{\partial \Omega} (a\cdot f_2(y)+g_2(y))^2dy+\lambda\|a\|_2^2
\end{equation}
and
\begin{equation}
F_n(a)=\frac{1}{n}\sum\limits_{i=1}^n l(a;x_i,y_i)+\lambda \|a\|_2^2
\end{equation}
respectively, where 
\begin{equation}
l(a;x_i,y_i)=|\Omega|(a\cdot f_1(x_i)+g_1(x_i))^2+|\partial \Omega|(a\cdot f_2(y_i)+g_2(y_i))^2
\end{equation}
and $i$-th components of $\mathbb{R}^m$-valued functions $f_1(x)$ and $f_2(y)$ are 
\begin{equation}
f_1^i(x)=-\|\omega_i\|_2^2\sigma^{''}(\omega_i \cdot x+b_i) +V(x)\sigma(\omega_i \cdot x+b_i)
\end{equation}
and $f_2^i(y)= \sigma(\omega_i\cdot y+b_i)$, respectively.

Note that the constructed random neural network $u$ has the form:
\begin{equation}
u(x)=\sum\limits_{i=1}^m a_i \sigma(\omega_i\cdot x+b_i),
\end{equation}
where $a=(a_1,\cdots, a_m)$, $\omega_1,\cdots, \omega_m \sim P_1$ and $b_1,\cdots, b_m \sim P_2$ (as defined in Theorem 1).

Since $P_1$ and $P_2$ are compactly supported in $\{\omega \in \mathbb{R}^d :\|\omega\|_1\leq M\}$ and $\{b\in \mathbb{R}: |b|\leq 2M\}$, we can deduce that for any $x\in \Omega$, $|f_1^i(x)|\lesssim M^2$ for any $i\in \{1,\cdots, m\}$ and then $\|f_1(x)\|_2 \lesssim \sqrt{m}M^2$.

Now, we are ready to present the generalization analysis for this method.

Let $\rho_0$ be a constant to be determined. Construct a sequence $\rho_k=2\rho_{k-1}$ for $k\in \mathbb{N}$ and then decompose $\mathcal{A}$ into a sequence of disjoint sets as follows.
\begin{equation}
\mathcal{A}_k:=\{a\in\mathcal{A}: \rho_{k-1} < F(a)-F(a^{*}) \leq \rho_k\},
\end{equation}
where $a^{*}= \mathop{\arg\min}_{a\in \mathcal{A}}  F(a)$. For brevity, we let $\rho_{-1}=0$ and $\mathcal{A}_0:=\{a\in\mathcal{A}: F(a)-F(a^{*}) \leq \rho_0\}$.

From $\|f_1(x)\|_2 \lesssim \sqrt{m}M^2$ and $\|f_2(y)\|_2 \lesssim 1$, we can deduce that
\begin{equation}
|F(a)-F(a^{*})|\lesssim (\lambda +m M^4) \|a-a^{*}\|_2 \lesssim \lambda +m M^4\lesssim m M^4,
\end{equation}
where the last inequality is due to the choice that $\lambda \in (0,1)$. This implies that the sequence is finite, specifically,
\begin{equation}
K:=\max k \leq C \log \left(\frac{m M^4}{\rho_0}\right),
\end{equation}
since $\rho_K=2^K\rho_0=\max_{a\in \mathcal{A}} F(a)-F(a^{*}) \leq \max_{a\in \mathcal{A}} F(a) \lesssim mM^4$. For any fixed constant $\delta \in (0,1)$, we set $\delta_k=\frac{\delta}{K+1}$ for $0\leq k \leq K$.

Let $s(a;x,y)=l(a;x,y)-l(a^{*};x,y)$, so that $F(a)-F_n(a)=(P-P_n)s(a;x,y)$, where $P_n$ is the empirical measure of $P$ and $P$ denotes the product measure of $Unif(\Omega)$ and $Unif(\partial \Omega)$, i.e., the uniform distribution on the interior and boundary, respectively.

By the fact that $F(\cdot)$ is $2\lambda$-strongly convex, we have $\lambda \|a-a^{*}\|_2^2\leq F(a)-F(a^{*})$. Thus, for any $a \in \mathcal{A}_k$, we have 
\begin{equation}
\|a-a^{*}\|_2^2\leq \frac{1}{\lambda} (F(a)-F(a^{*})) \leq \frac{\rho_k}{\lambda} .
\end{equation}

Due to the utilization of localization technique, we need to perform the generalization analysis in each subset $\mathcal{A}_k$. 

For any $a$ in $\mathcal{A}_k$ and $(x,y)\in \Omega\times \partial \Omega$, we know
\begin{equation}
\begin{aligned}
&|s(a;x,y)|\\
&= ||\Omega|(a\cdot f_1(x)+g_1(x))^2+|\partial\Omega|(a\cdot f_2(y)+g_2(y))^2-|\Omega|(a^{*}\cdot f_1(x)+g_1(x))^2-|\partial \Omega|(a^{*}\cdot f_2(y)+g_2(y))^2| \\
&\lesssim |(a-a^{*})\cdot f_1(x)| |a\cdot f_1(x)+g_1(x)+a^{*}\cdot f_1(x)+g_1(x) |  \\
&\quad+ |(a-a^{*})\cdot f_2(y)| |a\cdot f_2(y)+g_2(y)+a^{*}\cdot f_2(y)+g_2(y)      | \\
& \leq \|a-a^{*}\|_2 \|f_1(x)\|_2 (\|a\|_2 \|f_1(x)\|_2 +\|a^{*}\|_2 \|f_1(x)\|_2+2|g_1(x)|) + \|a-a^{*}\|_2 \|f_2(y)\|_2 (\|a\|_2 \|f_2(y)\|_2 +\|a^{*}\|_2 \|f_2(y)\|_2+2|g_2(y)|) \\
&\lesssim (\sqrt{\frac{\rho_k}{\lambda} } \sqrt{m}M^2)(\sqrt{m} M^2)\\
&= mM^4\sqrt{\frac{\rho_k}{\lambda} },
\end{aligned}
\end{equation}
where the last inequality following from that $\|a\|_2,\|a^{*}\|_2\lesssim 1, \|f_1(x)\|_2\lesssim  \sqrt{m}M^2$ and $\|f_2(y)\|_2\lesssim 1$. Note that here, for brevity, we have omitted the terms $|\Omega|$ and $|\partial \Omega|$.

With the upper bound for $|s(a;x,y)|$, applying McDiarmid’s inequality yields that with probability at least $1-\delta_k$, 
\begin{equation}
\sup\limits_{a\in \mathcal{A}_k} (P-P_n)s(a;x,y) \lesssim \mathbb{E}_{X,Y}\left[\sup\limits_{a\in \mathcal{A}_k} (P-P_n)s(a;x,y)\right]+ mM^4\sqrt{\frac{\rho_k}{\lambda} } \sqrt{\frac{\log(\frac{1}{\delta_k})}{n} }.
\end{equation}

By the technique of symmetrization, we have 
\begin{equation}
\begin{aligned}
\mathbb{E}_{X,Y}\left[\sup\limits_{a\in \mathcal{A}_k} (P-P_n)s(a;x,y)\right]&\leq 2\mathbb{E}_{X,Y,\epsilon} \left[\sup\limits_{a\in \mathcal{A}_k} \frac{1}{n}\sum\limits_{i=1}^n \epsilon_i (l(a;x_i,y_i)-l(a^{*};x_i,y_i))\right]\\
&\lesssim \mathbb{E}_{X,\epsilon}\left[\sup\limits_{a\in \mathcal{A}_k} \frac{1}{n}\sum\limits_{i=1}^n \epsilon_i ((a\cdot f_1(x_i)+g_1(x_i))^2 -(a^{*}\cdot f_1(x_i)+g_1(x_i))^2 )\right] \\
&\quad + \mathbb{E}_{Y,\epsilon}\left[\sup\limits_{a\in \mathcal{A}_k} \frac{1}{n}\sum\limits_{i=1}^n \epsilon_i ((a\cdot f_2(y_i)+g_2(y_i))^2 -(a^{*}\cdot f_2(y_i)+g_2(y_i))^2 )\right],
\end{aligned}
\end{equation}
where $X=(x_1,\cdots, x_n)$ and $Y=(y_1,\cdots, y_n)$ are samples from interior and boundary, respectively.

Then, from the contraction property of Rademacher complexity, we can deduce that
\begin{equation}
\begin{aligned}
&\mathbb{E}_{X,\epsilon}\left[\sup\limits_{a\in \mathcal{A}_k} \frac{1}{n}\sum\limits_{i=1}^n \epsilon_i ((a\cdot f_1(x_i)+g_1(x_i))^2 -(a^{*}\cdot f_1(x_i)+g_1(x_i))^2 )\right] \\
&\lesssim \sqrt{m}M^2 \mathbb{E}_{X,\epsilon}\left[\sup\limits_{a\in \mathcal{A}_k} \frac{1}{n}\sum\limits_{i=1}^n \epsilon_i (a-a^{*})\cdot f_1(x_i)\right] \\
&= \sqrt{m}M^2 \mathbb{E}_{X,\epsilon}\left[\sup\limits_{a\in \mathcal{A}_k} \frac{1}{n} (a-a^{*})\cdot \left(\sum\limits_{i=1}^n \epsilon_i f_1(x_i)\right)\right] \\
&\leq \sqrt{m}M^2 \mathbb{E}_{X,\epsilon}\left[\sup\limits_{a\in \mathcal{A}_k} \frac{1}{n} \|a-a^{*}\|_2 \left\| \sum\limits_{i=1}^n \epsilon_i f_1(x_i)\right\|_2\right] \\
&\leq \sqrt{m}M^2  \sqrt{\frac{\rho_k}{\lambda} }  \mathbb{E}_{\epsilon} \left[ \left\|\sum\limits_{i=1}^n \epsilon_i f_1(x_i)\right\|_2 \right] \\
&\leq \sqrt{m}M^2  \sqrt{\frac{\rho_k}{\lambda} }  \left( \mathbb{E}_{\epsilon}  \left\|\sum\limits_{i=1}^n \epsilon_i f_1(x_i)\right\|_2^2 \right)^{\frac{1}{2}} \\
&\lesssim \sqrt{m}M^2\frac{1}{n}  \sqrt{\frac{\rho_k}{\lambda} } \sqrt{n}\sqrt{m}M^2 \\
&=\frac{m M^4}{\sqrt{n}}  \sqrt{\frac{\rho_k}{\lambda} }.
\end{aligned}
\end{equation}

Similarly, we can derive the upper bound for the second term in (90). 

Thus, we have that with probability at least $1-\delta_k$,
\begin{equation}
\sup\limits_{a\in \mathcal{A}_k} (P-P_n)s(a;x,y) \lesssim \frac{mM^4}{\sqrt{n}} \sqrt{\frac{\rho_k}{\lambda} }+\frac{mM^4}{\sqrt{n}}\sqrt{\frac{\rho_k}{\lambda} } \sqrt{\log\left(\frac{1}{\delta_k}\right) }.
\end{equation}

Note that for any $a\in \mathcal{A}_k$, we have
\begin{equation}
\frac{\rho_k}{2}=\rho_{k-1}\leq F(a)-F(a^{*}),
\end{equation}
thus 
\begin{equation}
\rho_k \leq \max\{ \rho_0, 2(F(a)-F(a^{*}))\}.
\end{equation}

For $\delta_k$, it's obvious that 
\begin{equation}
\frac{1}{\delta_k}=\frac{K+1}{\delta}\leq C\frac{1}{\delta}\log\left(\frac{mM^4}{\rho_0}\right).
\end{equation}

Combining (92), (93), (95) with the fact that
\begin{equation}
(F(a)-F(a^{*}))-(F_n(a)-F_n(a^{*}))=(P-P_n)s(a;x,y)
\end{equation}
yields that
with probability at least $1-\delta_k$ for all $a\in \mathcal{A}_k$,
\begin{equation}
F(a)-F(a^{*}) \leq F_n(a)-F_n(a^{*}) +CmM^4\sqrt{\frac{\max\{ \rho_0, 2(F(a)-F(a^{*}))\}}{n\lambda}}\left( 1+\sqrt{\log\left(\frac{1}{\delta}\right)+\log\log \left(\frac{mM^4}{\rho_0}\right) }\right).
\end{equation}

Since $\sum_{k=0}^K \delta_k=\delta$, the above inequality holds with probability at least $1-\delta$ uniformly for all $a\in \mathcal{A}$.

By setting $\rho_0=1/n$, we can deduce that with probability at least $1-\delta$, for any $a\in \mathcal{A}$, we have 
\begin{equation}
F(a)-F(a^{*}) \leq \frac{\rho_0}{2}=\frac{1}{2n}
\end{equation}
or 
\begin{equation}
F(a)-F(a^{*}) \leq F_n(a)-F_n(a^{*}) +CmM^4\sqrt{\frac{2(F(a)-F(a^{*}))}{n\lambda}}\left( 1+\sqrt{\log\left(\frac{1}{\delta}\right)+\log\log \left(\frac{mM^4}{\rho_0}\right) }\right).
\end{equation}

From the form of $F_n(\cdot)$, we know that it is $2\lambda$-strongly convex, $CmM^4$-Lipschitz continuous and $CmM^4$-smooth. Thus, Theorem 3.10 in \cite{3} implies that the optimization error of projected gradient descent described in (31) with $\nu_t=\frac{1}{cmM^4}$ can be bounded as
\begin{equation}
\Delta_T:=F_n(a_T)-\min_{a\in\mathcal{A}} F_n(a)\lesssim \frac{(mM^4)^3}{\lambda^2}\exp{\left(-\frac{\lambda T}{CmM^4}\right)},
\end{equation}
where $a_T$ is obtained from (31) after $T$ iterations.

Therefore, combining (98), (99) and (100), we obtain that
\begin{equation}
\begin{aligned}
&F(a_T)-F(a^{*})\\
&\leq F_n(a_T)-F_n(a^{*})+CmM^4\sqrt{\frac{2(F(a_T)-F(a^{*}))}{n\lambda}}\left( 1+\sqrt{\log\left(\frac{1}{\delta}\right)+\log\log \left(\frac{mM^4}{\rho_0}\right) }\right) \\
&= F_n(a_T)-\min_{a\in\mathcal{A}} F_n(a)+\min_{a\in\mathcal{A}} F_n(a)-F_n(a^{*}) \\
&\quad+CmM^4\sqrt{\frac{2(F(a_T)-F(a^{*}))}{n\lambda}}\left( 1+\sqrt{\log\left(\frac{1}{\delta}\right)+\log\log \left(\frac{mM^4}{\rho_0}\right) }\right) \\
&\leq\Delta_T+CmM^4\sqrt{\frac{2(F(a_T)-F(a^{*}))}{n\lambda}}\left( 1+\sqrt{\log\left(\frac{1}{\delta}\right)+\log\log \left(\frac{mM^4}{\rho_0}\right) }\right) \\
&\leq \Delta_T+ \frac{F(a_T)-F(a^{*})}{2} + \frac{Cm^2M^8}{n\lambda}\left(1+\log\left(\frac{1}{\delta}\right)+\log\log \left(\frac{mM^4}{\rho_0}\right) \right),
\end{aligned}
\end{equation}
where the last inequality is from the basic inequality that $2\sqrt{ab}\leq a+b$ for $a,b\geq 0$.

Through a simple algebraic transformation for (101), we have
\begin{equation}
\begin{aligned}
F(a_T)-F(a^{*}) &\leq C\left(\frac{m^3 M^{12}}{\lambda^2}\exp{\left(-\frac{\lambda T}{CmM^4}\right)} + \frac{m^2M^8}{n\lambda}\left(1+\log\left(\frac{1}{\delta}\right)+\log\log \left(\frac{mM^4}{\rho_0}\right) \right)  \right).
\end{aligned}
\end{equation}

Integrating the tail term with respect to the samples from $\Omega$ and $\partial \Omega$ yields that
\begin{equation}
\mathbb{E}_{X,Y}[F(a_T)-F(a^{*})] \lesssim \frac{m^3 M^{12}}{\lambda^2}\exp{\left(-\frac{\lambda T}{CmM^4}\right)}+ \frac{m^2M^8}{n\lambda}\log\log \left(\frac{mM^4}{\rho_0}\right).
\end{equation}

Thus, it remains only to bound $F(a^{*})$. 

Note that since $u^{*} \in \mathcal{B}^s(\Omega) $ with $s\geq\frac{d+7}{2}$, Theorem 1 implies that there is a random neural network $u_m$ of the form
\begin{equation}
u_m(x)=\sum\limits_{i=1}^m \bar{a}_i \sigma(W_i\cdot x+B_i)
\end{equation}
such that
\begin{equation}
\mathbb{E}_{W,B} \left[  \|u_m(x)-u^{*}(x)\|_{H^2(\Omega)}^2\right] \lesssim \frac{M^{d+1}}{m}+\frac{1}{M^{d+1}},
\end{equation}
where $\bar{a}=(\bar{a}_1,\cdots, \bar{a}_m) \in \mathcal{A}$.

Therefore, from the definition of $a^{*}$, we have
\begin{equation}
\begin{aligned}
F(a^{*})&=\int_{\Omega} (a^{*} \cdot f_1(x)+g_1(x))^2dx+\int_{\partial\Omega} (a^{*} \cdot f_2(y)+g_2(y)))^2dy +\lambda\|a^{*}\|_2^2 \\
&\leq \int_{\Omega} (\bar{a} \cdot f_1(x)+g_1(x))^2dx+\int_{\partial\Omega} (\bar{a} \cdot f_2(y)+g_2(y)))^2dy +\lambda\|\bar{a}\|_2^2\\
&\lesssim \|u_m(x)-u^{*}(x)\|_{H^2(\Omega)}^2 +\lambda,
\end{aligned}
\end{equation}
which implies that 
\begin{equation}
\begin{aligned}
F(a^{*})&\lesssim \mathbb{E}_{W,B} \left[ \|u_m(x)-u^{*}(x)\|_{H^2(\Omega)}^2\right] +\lambda \\
&\lesssim \frac{M^{d+1}}{m}+\frac{1}{M^{d+1}}+\lambda.
\end{aligned}
\end{equation}

Taking $M$ such that $M^{2(d+1)}=m$, i.e., $M=m^{\frac{1}{2(d+1)}}$, from (103) and (107), we obtain that 
\begin{equation}
\begin{aligned}
&\mathbb{E}_{W,B}\mathbb{E}_{X,Y} \left[\int_{\Omega} (a_T \cdot f_1(x)+g_1(x))^2dx+\int_{\partial\Omega} (a_T \cdot f_2(y)+g_2(y)))^2dy\right] \\
&\leq\mathbb{E}_{W,B}\mathbb{E}_{X,Y} \left[ F(a_T)-F(a^{*})+F(a^{*}) \right] \\
&\lesssim \frac{m^3 M^{12}}{\lambda^2}\exp{\left(-\frac{\lambda T}{CmM^4}\right)}+ \frac{m^2M^8}{n\lambda}\log\log \left(\frac{mM^4}{\rho_0}\right) +\frac{M^{d+1}}{m}+\frac{1}{M^{d+1}}+\lambda \\
&= \frac{m^{3+\frac{6}{d+1}}}{\lambda^2} \exp{\left(-\frac{\lambda T}{Cm^{2+4/(d+1)}}\right)} + \frac{m^{2+\frac{4}{d+1}}}{\lambda n}\log\log(mn)+\frac{1}{\sqrt{m}}+\lambda.
\end{aligned}
\end{equation}

By choosing $\lambda$ and $m$ such that
\begin{equation}
\lambda=\frac{1}{\sqrt{m}}=\frac{m^{2+\frac{4}{d+1}}}{\lambda n},
\end{equation}
i.e., $\lambda=\frac{1}{\sqrt{m}}$ and $m=n^{\frac{d+1}{3d+7}}$, the final bound in (108) becomes 
\begin{equation}
\frac{m^{3+\frac{6}{d+1}}}{\lambda^2} \exp{\left(-\frac{\lambda T}{Cm^{1+2/(d+1)}}\right)} + \left( \frac{1}{n}\right)^{\frac{d+1}{2(3d+7)}}\log\log(n).
\end{equation} 
Thus, we can take $T$ such that the two terms in (110) are equal, i.e., 
\begin{equation}
T=\mathcal{O}\left( \frac{m^{1+2/(d+1)}}{\lambda} \log (n) \right)=\mathcal{O}(\sqrt{n}\log(n)).
\end{equation}
\end{proof}

\section{Proof of Theorem 4}
\begin{proof}
From the definition of $a_n$, we have $F_n(a_n)\leq F_n(0)$, thus $\lambda\|a_n\|_2^2 \leq C$, i.e., $\|a_n\|_2\leq \sqrt{\frac{C}{\lambda}}$.

Recall that 
\begin{equation}
Pl(a;x,y)=\int_{\Omega} (a \cdot f_1(x)+g_1(x))^2dx+\int_{\partial\Omega} (a \cdot f_2(y)+g_2(y)))^2dy.
\end{equation}

Therefore,
\begin{equation}
\begin{aligned}
\mathbb{E}_{X,Y}[Pl(a_n;x,y) ] &= \mathbb{E}_{X,Y}[(P-P_n)l(a_n;x,y) +P_nl(a_n;x,y) ] \\
&\leq \mathbb{E}_{X,Y} \left[\sup\limits_{\|a\|_2\leq \sqrt{\frac{C}{\lambda}}}(P-P_n)l(a;x,y) +P_nl(a_n;x,y)\right] \\
&\leq \mathbb{E}_{X,Y} \left[\sup\limits_{\|a\|_2\leq \sqrt{\frac{C}{\lambda}}}(P-P_n)l(a;x,y) +P_nl(a_n;x,y)+\lambda \|a_n\|_2^2\right] \\
&\leq \mathbb{E}_{X,Y} \left[\sup\limits_{\|a\|_2\leq \sqrt{\frac{C}{\lambda}}}(P-P_n)l(a;x,y) +P_nl(\bar{a};x,y)+\lambda \|\bar{a}\|_2^2\right] \\
&= \mathbb{E}_{X,Y} \left[\sup\limits_{\|a\|_2\leq \sqrt{\frac{C}{\lambda}}}(P-P_n)l(a;x,y)\right] +Pl(\bar{a};x,y)+\lambda \|\bar{a}\|_2^2\\
&= \mathbb{E}_{X,Y} \left[\sup\limits_{\|a\|_2\leq \sqrt{\frac{C}{\lambda}}}(P-P_n)l(a;x,y)\right]+ F(\bar{a}),
\end{aligned}
\end{equation}
where $\bar{a}$ is the output vector of the random neural network that provides the approximation result in Theorem 2 and the third inequality follows from the definition of $a_n$.

By the technique of symmetrization, we have
\begin{equation}
\begin{aligned}
&\mathbb{E}_{X,Y} \left[\sup\limits_{\|a\|_2\leq \sqrt{\frac{C}{\lambda}}}(P-P_n)l(a;x,y)\right] \\
&\lesssim 	\mathbb{E}_{X,\epsilon}\left[\sup\limits_{\|a\|_2\leq \sqrt{\frac{C}{\lambda}}} \frac{1}{n}\sum\limits_{i=1}^n\epsilon_i(a\cdot f_1(x_i)+g_1(x_i))^2\right]+	\mathbb{E}_{Y,\epsilon}\left[\sup\limits_{\|a\|_2\leq \sqrt{\frac{C}{\lambda}}} \frac{1}{n}\sum\limits_{i=1}^n\epsilon_i(a\cdot f_2(y_i)+g_2(y_i))^2\right].
\end{aligned}
\end{equation}

Note that the constructed random neural network $u$ is of the form:
\begin{equation}
u(x)=\sum\limits_{i=1}^m a_i\sigma(W_i\cdot x+B_i),
\end{equation}
where $a=(a_1,\cdots, a_m)$ and $W_1,\cdots, W_m\sim P_1$, $B_1,\cdots, B_m\sim P_2$ (as defined in Theorem 2).

Thus, $i$-component of $f_1(x)$ is 
\begin{equation}
f_1^i(x)=-\|W_i\|_2^2\sigma^{''}(W_i\cdot x+B_i)+V(x)\sigma(W_i\cdot x+B_i),
\end{equation}
which implies that for any $x\in \Omega$,
\begin{equation}
\|f_1(x)\|_2 \lesssim \sqrt{\sum\limits_{j=1}^m (\|W_j\|_2^4+1)}.
\end{equation}

Thus, for the first term in (114), we have
\begin{equation}
\begin{aligned}
&\mathbb{E}_{X} \mathbb{E}_{\epsilon}\left[\sup\limits_{\|a\|_2\leq \sqrt{\frac{C}{\lambda}}} \frac{1}{n}\sum\limits_{i=1}^n\epsilon_i(a\cdot f_1(x_i)+g_1(x_i))^2\right]  \\
&\lesssim \mathbb{E}_{X} \mathbb{E}_{\epsilon}\left[\sup\limits_{\|a\|_2\leq \sqrt{\frac{C}{\lambda}}} \frac{1}{n}\sum\limits_{i=1}^n\epsilon_i(a\cdot f_1(x_i)+g_1(x_i))\frac{1}{\sqrt{\lambda}} \sqrt{\sum\limits_{j=1}^m (\|W_j\|_2^4+1)} \right] \\ 
&=\frac{1}{n}\frac{1}{\sqrt{\lambda}}\sqrt{\sum\limits_{j=1}^m (\|W_j\|_2^4+1)}\mathbb{E}_{X} \mathbb{E}_{\epsilon}\left[\sup\limits_{\|a\|_2\leq \sqrt{\frac{C}{\lambda}}} \frac{1}{n} a\cdot \left(\sum\limits_{i=1}^n \epsilon_if_1(x_i)\right)  \right] \\
&\leq \frac{1}{n}\frac{1}{\sqrt{\lambda}}\sqrt{\sum\limits_{j=1}^m (\|W_j\|_2^4+1)}\mathbb{E}_{X} \mathbb{E}_{\epsilon} \left[\frac{1}{\sqrt{\lambda}} \left\|\sum\limits_{i=1}^n \epsilon_if_1(x_i) \right\|_2 \right]\\
&\leq \frac{1}{n}\frac{1}{\lambda}\sqrt{\sum\limits_{j=1}^m (\|W_j\|_2^4+1)}\mathbb{E}_{X} \left(\mathbb{E}_{\epsilon}  \left\|\sum\limits_{i=1}^n \epsilon_if_1(x_i) \right\|_2^2\right) ^{1/2}\\
&=\frac{1}{n}\frac{1}{\lambda}\sqrt{\sum\limits_{j=1}^m (\|W_j\|_2^4+1)}\mathbb{E}_{X} \left( \sum\limits_{i=1}^n \|f_1(x_i) \|_2^2\right) ^{1/2}\\
&\lesssim \frac{1}{n}\frac{1}{\lambda}\sqrt{\sum\limits_{j=1}^m (\|W_j\|_2^4+1)} \left(n \sum\limits_{j=1}^m (\|W_j\|_2^4+1)\right)^{1/2} \\
&=\frac{1}{\lambda\sqrt{n}} \sum\limits_{j=1}^m (\|W_j\|_2^4+1),
\end{aligned}
\end{equation}
where the first inequality follows from the contraction property of the Rademacher complexity and (117).

Similarly, for the second term in (114), we can deduce that
\begin{equation}
\mathbb{E}_{Y,\epsilon}\left[\sup\limits_{\|a\|_2\leq \sqrt{\frac{C}{\lambda}}} \frac{1}{n}\sum\limits_{i=1}^n\epsilon_i(a\cdot f_2(y_i)+g_2(y_i))^2\right] \lesssim \frac{m}{\lambda\sqrt{n}}.
\end{equation}

Since $W_1,\cdots, W_m\sim P_1$ and $P_1$ has density function $C_{\alpha}/(1+\|\omega\|_2)^{\alpha}$ with $\alpha>d+4$, we can deduce that $\mathbb{E}_{W} [\|W\|_2^4]\lesssim 1$. Combining this with (114), (118), (119), we have
\begin{equation}
 \mathbb{E}_{X,Y} \left[\sup\limits_{\|a\|_2\leq \sqrt{\frac{C}{\lambda}}}(P-P_n)l(a;x,y)\right] \lesssim \frac{m}{\lambda\sqrt{n}}.
\end{equation}

Thus, from (113), it remains only to bound $F(\bar{a})$. 

Recall that the constructed random neural network $u_m$ for approximating $u^{*}$ in Theorem 2 has the form 
\begin{equation}
u_m(x)=\sum\limits_{i=1}^m \bar{a}_i\sigma(W_i\cdot x+B_i),
\end{equation}
where 
\begin{equation}
\bar{a}_i= \frac{1}{m}\frac{1}{2\pi \hat{\sigma}(1) }\frac{|\hat{u^{*}}(W_i)|\cos(\theta(W_i)-B_i)}{p_1(W_i)p_2(B_i)} I_{\{|B_i|\leq Cd(1+\|W_i\|_2)\}}, 1\leq i \leq m.
\end{equation}

Thus for $F(\bar{a})$, we have
\begin{equation}
\begin{aligned}
\mathbb{E}_{W,B}[F(\bar{a})] &=\mathbb{E}_{W,B}\left[ \int_{\Omega} (\bar{a} \cdot f_1(x)+g_1(x))^2dx+\int_{\partial\Omega} (\bar{a} \cdot f_2(y)+g_2(y))^2dy\right]+\lambda \mathbb{E}_{W,B}[\|\bar{a}\|_2^2] \\
&\lesssim \mathbb{E}_{W,B}[\|u_m-u^{*}\|_{H^2(\Omega)}^2] +\lambda \mathbb{E}_{W,B}[\|\bar{a}\|_2^2]\\
&\lesssim \frac{1}{m}\frac{d^{\beta+1}}{C_{\alpha}} +\lambda \mathbb{E}_{W,B}[\|\bar{a}\|_2^2],
\end{aligned}
\end{equation}
where the first inequality follows from the Sobolev trace theorem and the second inequality is from Theorem 2.

Recall that the density functions of $P_1$ and $P_2$ are $C_{\alpha}/(1+\|\omega\|_2)^{\alpha}$ and $C_{\beta}/(1+|b|)^{\beta}$, respectively. Thus, combining with (122), we have that for any $i\in \{1,\cdots,m\}$,

\begin{equation}
\begin{aligned}
\mathbb{E}_{W,B}[\|\bar{a}_i\|_2^2]&\lesssim \frac{1}{m^2}\mathbb{E}_{W,B} \left[\frac{1}{C_{\alpha}^2C_{\beta}^2} |\hat{u^{*}}(W_i)|^2(1+\|W_i\|_2)^{2\alpha} (1+|B_i|)^{2\beta}I_{\{|B_i|\leq Cd(1+\|W_i\|_2)\}}\right]\\
&\lesssim \frac{1}{m^2}\mathbb{E}_{W,B} \left[\frac{d^{2\beta}}{C_{\alpha}^2C_{\beta}^2} |\hat{u^{*}}(W_i)|^2(1+\|W_i\|_2)^{2\alpha} (1+\|W_i\|_2)^{2\beta}\right]\\
&=\frac{1}{m^2}\mathbb{E}_{W} \left[\frac{d^{2\beta}}{C_{\alpha}^2C_{\beta}^2} |\hat{u^{*}}(W_i)|^2(1+\|W_i\|_2)^{2\alpha} (1+\|W_i\|_2)^{2\beta}\right]\\ 
&= \frac{1}{m^2}\frac{d^{2\beta}}{C_{\alpha}C_{\beta}^2} \int_{\mathbb{R}^d} (1+\|\omega\|_2)^{\alpha+2\beta}|\hat{u^{*}}(\omega)|^2d\omega\\
&\leq  \frac{1}{m^2}\frac{d^{2\beta}}{C_{\alpha}C_{\beta}^2} \int_{\mathbb{R}^d} (1+\|\omega\|_2)^{\alpha+\beta+5}|\hat{u^{*}}(\omega)|^2d\omega.
\end{aligned}
\end{equation}

Plugging (124) into (123) yields that
\begin{equation}
\mathbb{E}_{W,B}[F(\bar{a})] \lesssim \frac{1}{m}\frac{d^{\beta+1}}{C_{\alpha}} + \frac{\lambda}{m} \frac{d^{2\beta}}{C_{\alpha}}.
\end{equation}

Finally, combining (113), (120) and (125) leads to that
\begin{equation}
\mathbb{E}_{W,B}\mathbb{E}_{X,Y}\left[\int_{\Omega} (a_n \cdot f_1(x)+g_1(x))^2dx+\int_{\partial\Omega} (a_n \cdot f_2(y)+g_2(y))^2dy\right] \lesssim \frac{1}{m}\frac{d^{2\beta}}{C_{\alpha}}+ \frac{\lambda}{m}\frac{d^{2\beta}}{C_{\alpha}}+\frac{m}{\lambda\sqrt{n}}.
\end{equation}

\end{proof}

\end{document}